\documentclass[12pt]{amsart}
\usepackage{amssymb,amscd,amsthm,amsmath,color}
\usepackage{fullpage}
\usepackage{tikz-cd}
\usepackage{epsfig}
\usepackage{graphicx}
\usepackage{xypic}
\usepackage{url}
\usepackage{mathtools}
\usepackage{hyperref}
\usepackage{mathrsfs}
\usepackage{enumerate}
\usepackage{pigpen}

\numberwithin{equation}{section}

\newcommand{\OO}{\mathcal{O}}

\newcommand{\FF}{\mathcal{F}}

\newcommand{\Hom}{\textsf{Hom}}

\newcommand{\Spec}{\operatorname{Spec}}

\newcommand{\Pic}{\operatorname{Pic}}
\newcommand{\rk}{\operatorname{rk}}

\newcommand{\Eff}{\overline{\operatorname{Eff}}}

\newcommand{\Nef}{\operatorname{Nef}}

\newcommand{\Cl}{\operatorname{Cl}}

\newcommand{\Mor}{\operatorname{Mor}}

\newcommand{\po}{\ar@{}[dr]|{\text{\pigpenfont R}}}
\newcommand{\pb}{\ar@{}[dr]|{\text{\pigpenfont J}}}

\newtheorem{theorem}{Theorem}[section]
\newtheorem{lemma}[theorem]{Lemma}
\newtheorem{proposition}[theorem]{Proposition}
\newtheorem{corollary}[theorem]{Corollary}

\theoremstyle{definition}
\newtheorem{definition}[theorem]{Definition}

\newtheorem{conjecture}[theorem]{Conjecture}
\newtheorem{remark}[theorem]{Remark}

\newtheorem{question}[theorem]{Question}
\newtheorem{example}[theorem]{Example}
\newtheorem{construction}[theorem]{Construction}

\definecolor{ao(english)}{rgb}{0.0, 0.5, 0.0}

\input xy
\xyoption{all}

\makeatletter
\def\dasharrowfill@#1#2#3#4{%
        $\m@th
        \thickmuskip0mu
        \medmuskip\thickmuskip
        \thinmuskip\thickmuskip
        \relax
        #4#1\mkern2mu
        \xleaders\hbox{$#4\mkern2mu#2\mkern2mu$}\hfill
        \mkern2mu
        #3$%
}

\def\dashrightarrowfill@{\dasharrowfill@\relbar\relbar\rightarrow}

\providecommand*\xdashrightarrow[2][]{%
  \ext@arrow 0055{\dashrightarrowfill@}{#1}{#2}}
\makeatother

\begin{document}

\title{Globally F-regular varieties are rationally chain connected}
%\author[]{}

\author{Eric Jovinelly}
\address{Department of Mathematics \\
Brown University \\
Box 1917 \\
151 Thayer Street \\
Providence, RI, 02912}
\email{eric\_jovinelly@brown.edu}

\author{Brian Lehmann}
\address{Department of Mathematics \\
Boston College  \\
Chestnut Hill, MA \, \, 02467}
\email{lehmannb@bc.edu}

\author{Eric Riedl}
\address{Department of Mathematics \\
University of Notre Dame  \\
255 Hurley Hall \\
Notre Dame, IN 46556}
\email{eriedl@nd.edu}

\begin{abstract}
We prove that projective globally $F$-regular varieties are rationally chain connected.  We deduce this from a study of the positivity of their tangent bundles and their purely inseparable covers.  Although some globally $F$-regular varieties admit non-vanishing pluridifferential forms, we show that each admits a unique purely inseparable crepant cover that does not have this pathology.  These covers carry very free curves, which can be degenerated to chains of rational curves through general pairs of points.
\end{abstract}

\maketitle

\section{Introduction}

A fundamental theme of birational geometry is the close relationship between positivity of the tangent bundle and properties of rational curves on a projective variety $X$.  Following the innovative ideas of \cite{Mori79}, we can expect $X$ to carry many rational curves whenever $T_{X}$ is ``positive.''  \cite{KMM92, Campana92} established an influential and precise version of this statement:  smooth Fano varieties over an algebraically closed field of  characteristic $0$ are rationally connected. This was later extended to klt Fano varieties in characteristic 0 by \cite{HM07, Zhang06}.

In positive characteristic the situation is more subtle.  In this context, one typically works with rational chain connectedness due to the existence of inseparable families of curves.  For example, in characteristic $p$ the arguments of \cite{KMM92} show that smooth Fano varieties are rationally chain connected but do not settle rational connectedness.  It is a long-standing open problem to further our understanding of rational curves on smooth and singular Fano varieties in characteristic $p$.

We focus on projective globally $F$-regular varieties.  As the positive characteristic counterpart of varieties of klt Fano type in characteristic $0$, globally $F$-regular varieties have many special geometric features: they are Cohen-Macaulay with pseudo-rational singularities (\cite{Smith97}, \cite{Smith00}), satisfy vanishing theorems (\cite{Smith00}, \cite{SS10}), have only finitely many quasi-\'etale covers (\cite{CRST18}), and so on.
Even so, the qualitative behavior of rational curves on globally $F$-regular varieties has remained mysterious.
Tucker asked whether globally $F$-regular varieties are rationally chain connected (\cite{AIM13}); prior work establishes this for most threefolds and for particular classes of Fano varieties.

\subsection{Main results}
Our first main result establishes rational chain connectedness of projective globally $F$-regular varieties in all dimensions and characteristics.  The following theorem extends \cite[Theorem 0.1]{GLPSTZ15}, \cite[Theorem 1.5]{GNT19}, and \cite[Theorem 1.1]{OPTWX26} which address dimension $3$ and characteristic $p>5$.

\begin{theorem} \label{theo:maintheorem1}
Every projective globally $F$-regular variety is rationally chain connected.
\end{theorem}

Our approach to Theorem \ref{theo:maintheorem1} is different from prior works: we do not make use of the MMP or deform to characteristic $0$.  Instead, we establish a rich interplay between the following three features of a projective globally $F$-regular variety $X$:
\begin{itemize}
    \item the positivity of the tangent bundle of $X$;
    \item the set of crepant finite covers of $X$; and
    \item the existence of very free curves on $X$ (possibly of positive genus).
\end{itemize}
The behavior of globally $F$-regular varieties in these areas is not as simple as their characteristic 0 counterparts. For instance, there are globally $F$-regular varieties that have non-vanishing pluridifferential forms. (The most famous examples are certain degree $p$ covers of projective space constructed in \cite{Kollar95}; see Section \ref{sect:examples} for more details.)  However, our second main result shows that every projective globally $F$-regular variety with pluridifferential forms admits a canonical purely inseparable cover which does not have this pathology. 

\begin{theorem} \label{theo:canonicalcover}
    Let $X$ be a projective globally $F$-regular variety.  There is a unique purely inseparable crepant cover $f: Y \to X$ such that the smooth locus $V \subset Y$ satisfies $H^{0}(V,\Omega_{V}^{\otimes j}) = 0$ for every $j \geq 1$.  Furthermore, $Y$ is also globally $F$-regular.
\end{theorem}

Theorem \ref{theo:canonicalcover} suggests that globally $F$-regular varieties without pluridifferential forms play a distinguished role.  Our third main result confirms this: Such varieties share several important properties with klt Fano varieties in characteristic $0$.

\begin{theorem} \label{theo:introequiv}
    Let $X$ be a projective globally $F$-regular variety.  Then the following are equivalent:
    \begin{enumerate}
        \item $H^{0}(U,\Omega_{U}^{\otimes j})=0$ for every $j \geq 1$ where $U \subset X$ is the smooth locus.
        \item $X$ has no non-trivial purely inseparable crepant covers.
        \item $X$ admits a very free curve (possibly of positive genus) in its smooth locus.
        \item The only quotients of $T_{X}$ with numerically trivial first Chern class are torsion sheaves.
    \end{enumerate}
\end{theorem}

A famous question of Koll\'ar (\cite{Kollar}) asks whether every smooth Fano variety admits a very free rational curve.  Recently \cite{WW26} answered this question in the negative: there is a smooth Fano fourfold that admits no free rational curves at all.  In light of this advance, it is crucial to identify a class of Fano varieties which have the ``best behavior'' from the viewpoint of rational curves.  Theorem \ref{theo:introequiv} provides a natural candidate: these equivalent conditions are necessary for the existence of very free rational curves and we expect them to be sufficient.  (We further discuss this question in Section \ref{sect:rationalcurves}.)  The following definition identifies this class. 

\begin{definition}\label{def: superior cover}
    Let $X$ be a projective globally $F$-regular variety.
    \begin{enumerate}
        \item  We say that $X$ is superior if its smooth locus carries no pluridifferential forms.
        \item The superior cover of $X$ is the purely inseparable finite morphism $f: Y \to X$ identified by Theorem \ref{theo:canonicalcover}.
    \end{enumerate}
\end{definition}

To prove our main theorems, we establish a collection of novel structural results in the intertwined areas mentioned earlier.  We summarize these in the subsequent sections.

\subsection{Crepant covers}

Globally $F$-regular varieties have a special property: for every finite cover $f: Y \to X$ the divisor $K_{Y/X}$ is effective (see Lemma \ref{lemm:relcaneffstrong}).  The following definition identifies the ``boundary case'' which is of particular importance.

\begin{definition} \label{defi:crepantcover}
    Let $X$ be a normal variety.  A finite cover $f: Y \to X$ is a surjective finite morphism from a normal variety $Y$.  A crepant cover is a finite cover such that $K_{Y/X}$ is linearly equivalent to $0$.
\end{definition}

Suppose $X$ is a projective globally $F$-regular variety.  A separable crepant cover is the same thing as a quasi-\'etale map.  Thus the separable crepant covers are controlled by the \'etale fundamental group of the smooth locus of $X$.  In particular, \cite[Corollary 5.7]{CRST18} shows there can only be finitely many such covers.  (Our work provides a new, geometric proof of this result; see Theorem \ref{theo:finiteetalefund}.)

Our primary focus is purely inseparable crepant covers.  Previous work has focused on a special type: quasi-torsors under a finite infinitesimal group scheme.  Such covers have been studied in \cite{CR22}, \cite{KV25}, \cite{Vial26}; in particular \cite[Theorem A]{Vial26} shows that for a globally $F$-regular variety any tower of quasi-torsors will eventually stabilize.

We prove two structural properties of purely inseparable crepant covers that mirror the structure of quasi-\'etale covers.  The first is a finiteness result.

\begin{theorem} \label{theo:finiteness}
    Let $X$ be a projective globally $F$-regular variety.  Then $X$ admits only finitely many purely inseparable crepant covers up to isomorphism.    
\end{theorem}

The second is a factoring result.  By taking the universal cover of the smooth locus, we obtain a ``largest'' quasi-\'etale cover of a globally $F$-regular variety.  The superior cover plays a similar role for purely inseparable crepant covers.

\begin{theorem} \label{theo:maintheorem2}
    Let $X$ be a projective globally $F$-regular variety. The superior cover $f: Y \to X$ factors through every purely inseparable crepant cover of $X$.
\end{theorem}

\subsection{Positivity of the tangent bundle}

Traditionally mathematicians have studied the cotangent bundle of characteristic $p$ Fano varieties from the perspective of pluridifferential forms.  We will instead take a modern approach that focuses on slope positivity of the tangent bundle with respect to movable curves, following the ideas of \cite{CP11,Peternell12,CP19, Ou23}.

Recently \cite[Theorem 7.12]{CP19} and \cite[Theorem 1.4]{Ou23} have established a form of slope positivity that holds for all Fano and Calabi-Yau varieties in characteristic $0$.  Our next result proves the analogue for projective globally $F$-split varieties in characteristic $p$.  In low dimension, there are previous results in this direction by \cite{Kawakami21}.

\begin{theorem} \label{theo:maintheorem3}
     Let $X$ be a projective globally $F$-split variety.  Then every quotient of $T_{X}$ has pseudo-effective first Chern class.
\end{theorem}

The quotients of $T_{X}$ with numerically trivial first Chern class are particularly important.  Via Ekedahl's correspondence (\cite{Ekedahl87}) they control the height $1$ purely inseparable crepant covers of $X$.  The nicest geometry occurs when $T_{X}$ admits no such quotients at all -- in this case we say that $T_{X}$ has the positive quotient property (see Definition \ref{defi:pqp}).  This property is captured in Theorem \ref{theo:introequiv}.(4) and is the key technical tool for establishing the equivalent conditions in the theorem.

\subsection{Rational curves}
The idea of using slope positivity of the tangent bundle to study rational curves in characteristic $p$ has its roots in \cite{Miyaoka85}.  Similar techniques have been developed in characteristic $p$ by \cite{MM86, KSCT07, Tian15, Langer15, BM16, Langer21}. 
For our purpose, the most relevant prior work is \cite{JLR25}, which connects slope positivity of the tangent bundle to the existence of very free curves.  More precisely, \cite[Theorem 1.7]{JLR25} shows that if $T_{Y}$ satisfies the positive quotient property then $Y$ admits positive genus curves in its smooth locus that deform freely.  In particular, we may find such a curve passing through a general pair of points in $Y$ and deform it until it breaks into a union of geometrically rational curves that connect these points.  In this way we obtain chains of rational curves on (purely inseparable covers of) $X$ through pairs of general points.  This shows that globally $F$-regular varieties are rationally chain connected and establishes Theorem \ref{theo:maintheorem1}. 

\subsection{Summary of paper}
The logical flow of the paper is exactly backwards from the order of the results presented in the introduction.

We begin in Section \ref{sect:tangentpositivity} by establishing the slope positivity of the tangent bundle.  If $T_{X}$ admitted a negative quotient, then via Ekedahl's theory of $p$-closed foliations one could construct a purely inseparable cover $f: Y \to X$ such that $K_{Y/X}$ is not pseudo-effective.  However, Lemma \ref{lemm:relcanoneff} shows that no such cover can exist, thus proving Theorem \ref{theo:maintheorem3}.

In Section \ref{sect:crepantcovers} we turn to the study of purely inseparable crepant covers of globally $F$-regular varieties.  By Ekedahl's correspondence, this is roughly equivalent to studying the quotients of $T_{X}$ with numerically trivial first Chern class. By combining our slope positivity result with $F$-splitting techniques, we prove Theorem \ref{theo:finiteness} and (a version of) Theorem \ref{theo:maintheorem2}.

In Section \ref{sect:superior}, we connect slope positivity of the tangent bundle to the non-existence of pluridifferential forms using the results of \cite{JLR25}.  This allows us to reformulate the results of Section \ref{sect:crepantcovers} in terms of superior covers and thus to prove Theorem \ref{theo:canonicalcover} and Theorem \ref{theo:introequiv}.

In Section \ref{sect:rationalcurves} we prove Theorem \ref{theo:maintheorem1} by combining results of \cite{JLR25} with our other main theorems.  We also discuss several conjectures about stronger versions of rational connectedness for globally $F$-regular varieties.

Finally, Section \ref{sect:examples} is dedicated to examples, focusing on toric varieties and the construction of \cite{Kollar95}.

\subsection{Acknowledgements}

We thank Kevin Tucker for his generous help with this project.  We also thank Brendan Hassett, Eric Larson, Joaqu\'in Moraga, Karl Schwede, Zhiyu Tian, Isabel Vogt, Joe Waldron, Jakub Witaszek, and Zheng Xu for helpful conversations.

\subsection{AI usage}

This paper was written with input from Anthropic's Claude Opus 5.5 LLM and from OpenAI's ChatGPT 6 Sol LLM.  Specifically, these LLMs made two contributions that directly affected the theorems in the introduction:
\begin{itemize}
    \item In an earlier draft, all theorems in the introduction required a lower bound on the characteristic.  The characteristic assumption entered the argument at only one step, the proof of Theorem \ref{theo:weakpqp}, which used a generalization of an argument of \cite[Theorem 3.2]{Langer15}.  We then deduced Corollary \ref{coro:fsplitstronglyss} as a corollary.  ChatGPT pointed out that we should reverse the order: there is a short characteristic-free proof of Corollary \ref{coro:fsplitstronglyss} which can then be used to deduce Theorem \ref{theo:weakpqp}.  We implemented this change.
    \item An earlier draft did not include Theorem \ref{theo:finiteness}.  We were missing one key step: Theorem \ref{theo:finmanyheight1}.  The proof of this result was communicated to us by ChatGPT and is based on the work of \cite{AE11,BS11}; we did not contribute to the ideas.
\end{itemize}
LLMs also made several additional contributions:
\begin{itemize}
    \item  In an earlier draft, our argument for Proposition \ref{prop:constructingpqpcover} combined $F$-splitting techniques with ideas from \cite{Shen10}.  Jakub Witaszek suggested this might be related to $F$-signatures, and then ChatGPT separately suggested we simplify the proof by using \cite{CR22}.  We implemented this suggestion. 
    \item ChatGPT provided a simple argument for Proposition \ref{prop:foliationcoverfsplitting}.
    \item Claude explained the relevance of \cite{Sun08} for analyzing purely inseparable covers of SRC varieties; we adapted this approach when proving Proposition \ref{prop:allarePICC}. 
    \item We used Claude to compute Example \ref{exam:kollar}.
    \item AIs proofread the paper, identifying an important correction to Construction \ref{cons:toricsuperior} as well as some typos.
\end{itemize}

\section{Background} \label{sect:background}

Throughout our ground field $\mathbb{K}$ will be algebraically closed of prime characteristic $p$.  A variety denotes an irreducible reduced separated scheme of finite type over the ground field.  An open subset of a variety is big if its complement has codimension at least $2$.

Given a $\mathbb{K}$-scheme $X$, we denote by $F_{abs}: X \to X$ the absolute Frobenius.  We denote by $F_{rel}: X \to X^{(1)}$ the relative Frobenius over $\Spec(\mathbb{K})$, i.e.~the $\mathbb{K}$-morphism in the diagram
\begin{equation*} %\label{eq:frob}
\xymatrix{
X \ar@/^2pc/[rr]^{F_{abs}} \ar[r]_{F_{rel}} \ar[dr] & X^{(1)} \pb \ar[r] \ar[d] & X \ar[d] \\
& \Spec(\mathbb{K}) \ar[r]_{F_{abs}} & \Spec(\mathbb{K})
}
\end{equation*}
By repeated applications of $F_{rel}$, we get a sequence of $\mathbb{K}$-morphisms which we denote by
\begin{equation*}
\ldots X^{(-2)} \xrightarrow{F_{rel}}  X^{(-1)} \xrightarrow{F_{rel}} X \xrightarrow{F_{rel}} X^{(1)} \xrightarrow{F_{rel}} X^{(2)} \xrightarrow{F_{rel}} \ldots
\end{equation*}

Suppose $X$ is a normal variety.  For a torsion-free sheaf $\mathcal{F}$ on $X$, we use $\widehat{\mathcal{F}}$ to denote the double dual $\mathcal{F}^{\vee \vee}$.  Given torsion-free sheaves $\mathcal{E},\mathcal{F}$ on $X$, we denote by $\mathcal{E} \widehat{\otimes} \mathcal{F}$ the double dual $(\mathcal{E} \otimes \mathcal{F})^{\vee \vee}$ and by $\widehat{\wedge}^{k} \mathcal{E}$ the double dual $(\wedge^{k} \mathcal{E})^{\vee \vee}$. Given a morphism $f: Y \to X$, we write $f^{\widehat{*}}\mathcal{E}$ for $(f^{*}\mathcal{E})^{\vee \vee}$.  Note that if $f: Y \to X$ is a finite surjective morphism of normal varieties then $f_{*}$ preserves reflexivity of coherent sheaves.

\subsection{Divisors}
Let $X$ be a normal projective variety.  We denote the Picard group by $\Pic(X)$ and the class group by $\Cl(X)$.

For any torsion-free coherent sheaf $\mathcal{F}$ on $X$, its first Chern class $c_{1}(\mathcal{F})$ is the element of $\Cl(X)$ obtained by choosing a big smooth open subset $U \subset X$ such that $\mathcal{F}|_{U}$ is locally free and then extending the line bundle $\bigwedge^{\rk(\mathcal{F})} \mathcal{F}|_{U}$ to a reflexive rank $1$ sheaf on $X$.  (We will frequently engage in a mild abuse of notation by letting $c_{1}(\mathcal{F})$ denote any Weil divisor representing this linear equivalence class.) Note that this construction is insensitive to the behavior of $\mathcal{F}$ along codimension $\geq 2$ closed subsets; for example, we have $c_{1}(\mathcal{F}) = c_{1}(\widehat{\mathcal{F}})$.

\begin{remark}\label{rem: pullback for finite covers}
    Suppose $f: Y \to X$ is a finite surjective morphism of normal varieties.  Then we can define a pullback map $f^{*}: \Cl(X) \to \Cl(Y)$ in the following way.  Since $f$ is finite and $Y$ and $X$ are normal, the $f$-image of the singular locus of $Y$ has codimension $\geq 2$ and the $f$-preimage of the singular locus of $X$ has codimension $\geq 2$.  Thus we can choose big smooth open subsets $V \subset Y$ and $U \subset X$ such that $f^{-1}(U) = V$.  We then define the pullback of a Weil divisor $D$ on $X$ to be the extension to $Y$ of $f^{*}D|_{U}$.  In particular, we can define the relative canonical divisor $K_{Y/X} = K_{Y} - f^{*}K_{X}$ as a Weil divisor (up to linear equivalence).
\end{remark}

\subsection{Numerical equivalence}

Suppose $X$ is a normal projective variety of dimension $n$.
We let $N_{n-1}(X)_{\mathbb{Z}}$ denote the additive group of Weil divisors up to numerical equivalence (that is, up to intersecting against homogeneous polynomials of weighted degree $(n-1)$ in Chern classes of vector bundles on $X$; see \cite{FL17}).  We let $N^{n-1}(X)_{\mathbb{Z}}$ denote the additive group of homogeneous polynomials in Chern classes of vector bundles of weighted degree $(n-1)$ up to numerical equivalence.  There is a natural intersection pairing $N^{n-1}(X)_{\mathbb{Z}} \times N_{n-1}(X)_{\mathbb{Z}} \to \mathbb{Z}$.  By tensoring these groups with $\mathbb{R}$ we get dual finite-dimensional vector spaces $N_{n-1}(X)_{\mathbb{R}}$ and $N^{n-1}(X)_{\mathbb{R}}$.

Given a generically finite morphism $f: X \to Y$ of normal projective varieties of dimension $n$, we get maps $f_{*}: N_{n-1}(X)_{\mathbb{Z}} \to N_{n-1}(Y)_{\mathbb{Z}}$ induced by pushforward of cycles and maps $f^{*}: N^{n-1}(Y)_{\mathbb{Z}} \to N^{n-1}(X)_{\mathbb{Z}}$ induced by pullback of vector bundles.

The pseudo-effective cone of Weil divisors $\Eff_{n-1}(X) \subset N_{n-1}(X)_{\mathbb{R}}$ is the closure of the cone generated by classes of effective Weil divisors.  This cone is full-dimensional, closed, convex, and pointed by \cite[Theorem 1.4.(1)]{FL17}.  Its dual is the nef cone of curves $\Nef^{n-1}(X) \subset N^{n-1}(X)_{\mathbb{R}}$.  By duality $\Nef^{n-1}(X)$ is also full-dimensional, closed, convex, and pointed.

\begin{remark}
    When $X$ is $\mathbb{Q}$-factorial, the spaces $N_{n-1}(X)_{\mathbb{R}}$ and $N^{n-1}(X)_{\mathbb{R}}$ and the cones they contain coincide with the (more familiar) numerical spaces constructed using the duality between Cartier divisors and effective $1$-cycles.
\end{remark}

\subsection{Positive curves}
The following definition identifies a higher-genus analogue of a very free rational curve. 

\begin{definition}\label{def: very free curve}
    Let $X$ be a projective variety and let $C$ be a smooth projective curve.  We say that a morphism $s: C \to X$ is very free if $s(C)$ lies in the smooth locus of $X$ and every non-zero torsion-free quotient of $s^{*}T_{X}$ has slope at least $2g(C)+1$.
\end{definition}

In this situation $s^{*}T_{X}$ is ample by \cite[Corollary 4.5]{Zhao17}.  \cite[Proposition 2.7]{JLR25} verifies that if we fix a closed subset $Z \subset X$ of codimension $\geq 2$ then a general deformation of $s$ will avoid $Z$.

\subsection{Frobenius split varieties}

We quickly recall two basic types of Frobenius split varieties following \cite{SS10}.  The first is the characteristic $p$ analogue of a variety of klt Fano type and the second is the analogue of a variety of lc Calabi-Yau type. 

\begin{definition}
    Let $X$ be a normal variety.  We say that:
    \begin{itemize}
    \item $X$ is globally $F$-regular if, for every effective divisor $D$, there is an $e>0$ such that there is a splitting of the map
    \begin{equation*}
        \mathcal{O}_{X} \to F^{e}_{abs,*}\mathcal{O}_{X}(D).
    \end{equation*}
    \item $X$ is globally $F$-split if for some $e>0$ (equivalently, all $e>0$) the map $\mathcal{O}_{X} \to F^{e}_{abs,*}\mathcal{O}_{X}$ splits.
    \end{itemize}
\end{definition}
There is an implication
\begin{equation*}
    X \textrm{  globally $F$-regular} \implies X \textrm{ globally $F$-split.}
\end{equation*}

\begin{lemma} \label{lemm:frobeniustwistandsplitting}
    Let $X$ be a normal projective variety and consider its Frobenius twists $X^{(i)}$, $i \in \mathbb{Z}$.
    \begin{enumerate}
        \item If $X$ is globally $F$-split, so is $X^{(i)}$ for every $i$.
        \item If $X$ is globally $F$-regular, so is $X^{(i)}$ for every $i$.
    \end{enumerate}
\end{lemma}

\begin{proof}
    (1) Since we are working over an algebraically closed field $\mathbb{K}$, each $X^{(i)}$ is abstractly isomorphic to $X$.  As the globally $F$-split condition does not rely on the $\mathbb{K}$-structure, the conclusion holds.

    (2) The argument is the same.
\end{proof}

\subsection{p-closed foliations}
Let $X$ be a normal variety.  A foliation $\mathcal{F} \subset T_{X}$ is a saturated subsheaf that is closed under the Lie bracket. A foliation $\mathcal{F}$ is said to be $p$-closed if it is closed under $p$-powers, that is, for every smooth open $U \subset X$ and for every derivation $D: \mathcal{O}_{U} \to \mathcal{O}_{U}$ in $\mathcal{F}(U)$ the $p$-fold composition $D^{p}$ remains in $\mathcal{F}(U)$.  

\begin{lemma}[{\cite[Proposition 4.2]{Ekedahl87}}] \label{lemm:ekedahl}
    Let $X$ be a normal variety and let $\mathcal{F} \subset T_{X}$ be a saturated subsheaf.
    \begin{enumerate}
        \item If $\Hom( \widehat{\wedge}^2 \mathcal{F} , \widehat{T_X / \mathcal{F}}) = 0$  then $\mathcal{F}$ is a foliation.
        \item If $\mathcal{F}$ is a foliation and $\Hom({F_{abs}^{\widehat{*}}\mathcal{F}}, \widehat{T_X / \mathcal{F}}) = 0$ then $\mathcal{F}$ is $p$-closed.
    \end{enumerate}
\end{lemma}

%%%%%KEEP IN TEX FILE
%\begin{proof}
%    (1) Let $U \subset X$ be a smooth open subset such that $X \backslash U$ has codimension $2$ and $\mathcal{F}|_{U}$ and $ (T_{X}/\mathcal{F})|_{U}$ are both locally free.  Our assumption implies that $\Hom( {\wedge}^2 \mathcal{F}|_{U}, (T_X / \mathcal{F})|_{U}) = 0$.  By \cite[Proposition 4.2.(i)]{Ekedahl87} we conclude that $\mathcal{F}|_{U}$ is closed under the Lie bracket. Since $T_{X}$ is reflexive, $\mathcal{F}$ is reflexive (being a saturated subsheaf of a reflexive sheaf), and $X \backslash U$ has codimension $\geq 2$, membership in $\mathcal{F}$ can be tested after restricting to $U$, showing that $\mathcal{F}$ is also closed under the Lie bracket.

%    (2) It is a similar argument using \cite[Proposition 4.2.(ii)]{Ekedahl87} in place of \cite[Proposition 4.2.(i)]{Ekedahl87}.
%\end{proof}

By \cite{Ekedahl87}, every $p$-closed foliation induces a height $1$ purely inseparable finite morphism $f: X \to Y$ to a normal variety $Y$.  Choose big smooth open subsets $U \subset X$ and $V \subset Y$ such that $f^{-1}(V) = U$ and $\mathcal{F}|_{U}$ is locally free.  In this situation, \cite[Corollary 3.4]{Ekedahl87} yields an exact sequence of locally free sheaves
\begin{equation} \label{eq:tangentondomain}
    0 \to \FF|_{U} \to T_U \to f^*T_{V} \to F_{abs}^*\FF|_{U} \to 0  .
\end{equation}
Since $f$ is height $1$, we also obtain a morphism $g: Y^{(-1)} \to X$.  We will normally focus on this latter map which we call the purely inseparable cover induced by $\mathcal{F}$.  On $V^{(-1)}$ we have an exact sequence
\begin{equation} \label{eq:twistedtangent}
    0 \to \mathcal{K} \to T_{V^{(-1)}} \to g^*T_{U} \to g^{*}(T_{X}/\mathcal{F})|_{U} \to 0
\end{equation}
where $\mathcal{K}$ is the kernel of the natural map $T_{V^{(-1)}} \to g^{*}\mathcal{F}|_{U}$ and satisfies $F_{abs}^{*}\mathcal{K} \cong g^{*}(T_{X}/\mathcal{F})|_{U}$.  By taking reflexive closures and then taking first Chern classes, we obtain:

\begin{lemma} \label{lemm:relCanonicalHeight1}
    Let $\FF$ be a $p$-closed foliation on a normal projective variety $X$. Then the induced purely inseparable cover $g: Y^{(-1)} \to X$ satisfies $pK_{Y^{(-1)}/X} = (p-1)g^{*}c_1(T_X/\FF)$.
\end{lemma}

%%%%%%KEEP IN TEX FILE
%\begin{proof}
%    Let $f: X \to Y$ denote the quotient by $\mathcal{F}$.  Choose smooth open subsets $U \subset X$, $V \subset Y$ such that $f^{-1}(V) = U$ and the complements have codimension $\geq 2$.  By Equation \eqref{eq:tangentondomain} we have an exact sequence 
%    \[ 0 \to \FF|_{U} \to T_U \to g^*T_{V} \to F_{abs}^*\FF|_{U} \to 0  .\]
%    Pulling back to $V^{(-1)}$, we get
%    \[ 0 \to g^* \FF|_{U} \to g^* T_U \to F_{rel}^* T_{V} \to g^* F_{abs}^* \FF|_{U} \to 0 , \]
%    Taking reflexive closures and then taking first Chern classes, we obtain an equality of Weil divisors
%    \[ - g^{*}c_1(\FF) -g^*K_X  + pK_Y^{(-1)} + p g^{*}c_1(\FF) = 0 \]
%    which we rearrange to get $pK_{Y^{(-1)}/X} = (p-1)g^{*}c_1(T_X / \FF)$.
%\end{proof}

\subsection{Slope stability}
\cite{CP11} develops the theory of semistability with respect to nef curve classes on smooth projective varieties.  This is extended to $\mathbb{Q}$-factorial normal projective varieties by \cite{GKP14}, \cite{GKP16} and to arbitrary normal projective varieties by \cite{JLR26}.  We briefly summarize the main features of this theory.

Let $X$ be a normal projective variety of dimension $n$ and let $\alpha \in \Nef^{n-1}(X)$.  For a non-zero torsion-free coherent sheaf $\mathcal{E}$ on $X$, its $\alpha$-slope is
\begin{equation*}
    \mu_{\alpha}(\mathcal{E}) = \frac{c_{1}(\mathcal{E}) \cdot \alpha}{\rk(\mathcal{E})}.
\end{equation*}
We omit the subscript $\alpha$ whenever the curve class is understood from the context.

We say that $\mathcal{E}$ is semistable if for every non-zero proper subsheaf $\mathcal{F} \subset \mathcal{E}$ we have $\mu(\mathcal{F}) \leq \mu(\mathcal{E})$.  We say that $\mathcal{E}$ is strongly semistable if $F^{e\widehat{*}}_{abs}\mathcal{E}$ is semistable for every positive integer $e$.  Note that both notions inherit insensitivity to changes in codimension $\geq 2$ via the analogous property of first Chern classes.

Every non-zero torsion-free coherent sheaf $\mathcal{E}$ admits a Harder-Narasimhan filtration, i.e.~a sequence of saturated subsheaves
\begin{equation*}
    0 = \mathcal{F}_{0} \subset \mathcal{F}_{1} \subset \ldots \subset \mathcal{F}_{t} = \mathcal{E}
\end{equation*}
such that each successive quotient $\mathcal{F}_{i}/\mathcal{F}_{i-1}$ is semistable and the slopes of the successive quotients are strictly decreasing.  We define:
\begin{itemize}
    \item $\mu^{max}(\mathcal{E})$ to be the maximal slope of any subsheaf, i.e.~$\mu^{max}(\mathcal{E}) = \mu(\mathcal{F}_{1})$;
    \item $\mu^{min}(\mathcal{E})$ to be the minimal slope of any torsion-free quotient, i.e.~$\mu^{min}(\mathcal{E}) = \mu(\mathcal{E}/\mathcal{F}_{t-1})$;
    \item the maximal positive subsheaf to be the last term $\mathcal{F}_{i}$ satisfying $\mu^{min}(\mathcal{F}_{i}) > 0$ (or the zero sheaf if no such term exists).
\end{itemize}

The basic properties of slope stability in the setting of arbitrary normal projective varieties equipped with nef curve classes are essentially no different from the setting of smooth projective varieties equipped with ample polarizations; see \cite{JLR26} for details.

\section{Foliations on globally $F$-split varieties} \label{sect:tangentpositivity}

Our first goal is to understand the positivity properties of the tangent bundle $T_X$ when $X$ is a projective globally $F$-split variety.  Our main results are motivated by the work of \cite{CP19} and \cite{Ou23} which study the analogous question in characteristic $0$ for lc pairs $(X,\Delta)$ with $-(K_{X}+\Delta)$ nef.  Some of the techniques used in this section were pioneered by \cite{Miyaoka85}.

\subsection{Positivity of the relative canonical divisor}

One ingredient is a fundamental property of purely inseparable covers of a globally $F$-split variety: the relative canonical divisor must be effective.

\begin{lemma} \label{lemm:relcanoneff}
Let $X$ be a globally $F$-split variety.  For any purely inseparable finite cover $f: Y \to X$ the divisor $K_{Y/X}$ is linearly equivalent to an effective Weil divisor.
\end{lemma}

\begin{proof}
Choose $e$ sufficiently large that $F_{rel}^{e}$ factors through $f$.  Thus we have a factoring
\begin{equation*}
\mathcal{O}_{X} \to f_{*}\mathcal{O}_{Y} \to F^{e}_{abs,*}\mathcal{O}_{X}
\end{equation*}
Since the outer map splits, $\mathcal{O}_{X} \to f_{*}\mathcal{O}_{Y}$ also splits.

Choose big smooth open subsets $V \subset Y$ and $U \subset X$ such that $f^{-1}(U) = V$. By tensoring a splitting of $\mathcal{O}_{U} \to f_{*}\mathcal{O}_{V}$ by $\omega_{U}$, we obtain a non-zero homomorphism $f_{*}(f^{*}\omega_{U}) \cong f_{*}\mathcal{O}_{V} \otimes \omega_{U} \to \omega_{U}$.  By Grothendieck duality,
\begin{equation*}
\mathcal{H}om(f_{*}(f^{*}\omega_{U}),\omega_{U}) \cong f_{*}\mathcal{H}om(f^{*}\omega_{U},\omega_{V}).
\end{equation*}
Thus we have a non-zero global section on the right hand side.  This shows that $K_{Y/X}$ is effective.
\end{proof}

There is a stronger version of Lemma \ref{lemm:relcanoneff} for globally $F$-regular pairs.

\begin{lemma}[{\cite[Lemma 6.14]{BMPSTWW23}}] \label{lemm:relcaneffstrong}
    Let $(X,\Delta)$ be a globally $F$-regular pair (as defined in \cite{SS10}).  For any finite cover $f: Y \to X$ the divisor $K_{Y/X} - \lfloor f^{*}\Delta \rfloor$ is linearly equivalent to an effective Weil divisor.
\end{lemma}

\begin{proof}
    \cite[Lemma 6.14]{BMPSTWW23} shows that $(X,\Delta)$ is globally +-regular, implying that $\mathcal{O}_{X} \to f_{*}\mathcal{O}_{Y}(\lfloor f^{*}\Delta \rfloor)$ splits.  The conclusion follows by applying Grothendieck duality as in the proof of Lemma \ref{lemm:relcanoneff}.
\end{proof}

\subsection{Strong semistability}

The other key input is a strong semistability result.

\begin{theorem} \label{theo:strongsemistability}
    Let $X$ be a normal projective variety and fix a nef curve class $\alpha$.  Suppose that every height $1$ purely inseparable cover $f: Y \to X$ satisfies $K_{Y/X} \cdot f^{*}\alpha \geq 0$.  Then every semistable torsion-free sheaf $\mathcal{E}$ on $X$ is strongly semistable.
\end{theorem}

\begin{proof}
    The proof is by induction on the rank $r$ of $\mathcal{E}$.  The base case $r=1$ is trivial.

    For the induction step, it suffices to show that $F_{abs}^{\widehat{*}}\mathcal{E}$ is semistable.  Suppose for a contradiction that this is not the case and let $\mathcal{G}$ denote the maximal destabilizing subsheaf.  If $\mathcal{G}$  descended under $F_{abs}$ this sheaf would destabilize $\mathcal{E}$, a contradiction.  Thus \cite[Theorem 5.1]{Katz70} shows that on the smooth locus $U \subset X$ the Cartier connection $\nabla: F_{abs}^{\widehat{*}}\mathcal{E}|_{U} \to \Omega_{U} \otimes F_{abs}^{\widehat{*}}\mathcal{E}|_{U}$ cannot map $\mathcal{G}|_{U}$ to $\Omega_{U} \otimes \mathcal{G}|_{U}$.  Thus the induced $\mathcal{O}_{U}$-linear map
    \begin{equation*}
        \mathcal{G}|_{U} \to \Omega_{U} \otimes F^{*}_{abs}\mathcal{E}/\mathcal{G}|_{U}
    \end{equation*}
    is non-trivial.  This corresponds to a non-zero morphism $\phi: T_{X} \to \widehat{\mathcal{H}om}(\mathcal{G},F^{\widehat{*}}_{abs}\mathcal{E}/\mathcal{G})$. 

    Let $\mathcal{F}$ denote the kernel of $\phi$.  We claim that $\mathcal{F}$ is a $p$-closed foliation.  Indeed, on any open subset $V \subset U$ we see that $\mathcal{F}|_{V}$ consists of the vector fields $D$ such that $\nabla_{D} (\mathcal{G}|_{V}) \subset \mathcal{G}|_{V}$.  Since $\nabla_{[D_{1},D_{2}]} = [\nabla_{D_{1}},\nabla_{D_{2}}]$, the Lie bracket preserves $\mathcal{F}|_{V}$.  Similarly, since $\nabla_{D^{[p]}} = \nabla_{D}^{p}$ we see that $\mathcal{F}|_{V}$ is $p$-closed.

    Let $f: Y \to X$ be the height $1$ purely inseparable cover corresponding to $\mathcal{F}$ and let $\mathcal{Q}$ denote $T_{X}/\mathcal{F}$.  Lemma \ref{lemm:relCanonicalHeight1} shows that
    \begin{equation*}
        K_{Y/X} \cdot f^{*}\alpha = \tfrac{\deg(f)(p-1)}{p} c_{1}(\mathcal{Q}) \cdot \alpha.
    \end{equation*}
    Note that $\mu(\mathcal{Q}) \leq \mu^{max}(\widehat{\mathcal{H}om}(\mathcal{G},F^{\widehat{*}}_{abs}\mathcal{E}/\mathcal{G}))$.  Furthermore, our induction assumption combined with \cite[Corollary 5.18]{JLR26} implies
    \begin{equation*}
        \mu^{max}(\widehat{\mathcal{H}om}(\mathcal{G},F^{\widehat{*}}_{abs}\mathcal{E}/\mathcal{G})) = \mu^{max}(F^{\widehat{*}}_{abs}\mathcal{E}/\mathcal{G}) -\mu^{min}(\mathcal{G}) < 0.
    \end{equation*}
    We conclude that $c_{1}(\mathcal{Q}) \cdot \alpha < 0$, contradicting our assumption that $K_{Y/X} \cdot f^{*}\alpha \geq 0$.
\end{proof}

By combining with Lemma \ref{lemm:relcanoneff} we immediately conclude:

\begin{corollary} \label{coro:fsplitstronglyss}
    Suppose $X$ is a projective globally $F$-split variety.  For any nef curve class $\alpha$, every semistable sheaf is strongly semistable.
\end{corollary}

\begin{corollary} \label{coro:pclosedcriterion}
    Suppose $X$ is a projective globally $F$-split variety and let $\alpha$ be a nef curve class.  Suppose that $\mathcal{F}$ is a term in the Harder-Narasimhan filtration of $T_{X}$ satisfying $\mu^{min}(\mathcal{F}) \geq 0$.  Then $\mathcal{F}$ is a $p$-closed foliation.
\end{corollary}

\begin{proof}
    Corollary \ref{coro:fsplitstronglyss} shows that $\mu^{min}(F_{abs}^{\widehat{*}}\mathcal{F}) = p \mu^{min}(\mathcal{F})$ and that $\mu^{min}(\widehat{\wedge}^{2}\mathcal{F}) \geq 2 \mu^{min}(\mathcal{F})$.  Thus both quantities are larger than $\mu^{max}(T_{X}/\mathcal{F})$.  We conclude by Lemma \ref{lemm:ekedahl}.
\end{proof}

\subsection{Positivity of the tangent bundle}

Our main result in this section shows that for a projective globally $F$-split variety $X$ the tangent bundle $T_{X}$ can have no ``negative'' quotients.  The following theorem is an analogue of the characteristic $0$ result \cite[Theorem 1.4]{Ou23}.

\begin{theorem} \label{theo:weakpqp}
    Suppose $X$ is a projective globally $F$-split variety.  Then for any non-zero quotient $T_X \to Q$ the first Chern class $c_1(Q)$ is pseudo-effective. 
\end{theorem}

\begin{proof}
    \cite[Theorem 4.3]{SS10} shows that $-K_{X}$ is $\mathbb{Q}$-linearly equivalent to an effective $\mathbb{Q}$-Weil divisor.  In particular $\mu(T_{X}) \geq 0$.
    
    Suppose there were a non-zero quotient $T_X \to Q$ with $c_1(Q) \cdot \alpha < 0$ for some nef class $\alpha \in N^{n-1}(X)_{\mathbb{R}}$.  Let $\mathcal{F}_i$ be the last term in the Harder-Narasimhan filtration of $T_X$ such that $\mu^{min}(\mathcal{F}_{i}) \geq 0$.  Corollary \ref{coro:pclosedcriterion} implies $\mathcal{F}_i$ is a $p$-closed foliation.  Hence, we can construct a height one purely inseparable cover $g:Z \to X$ corresponding to $\mathcal{F}_{i}$. By Lemma \ref{lemm:relCanonicalHeight1}, we see that $K_{Z/X}$ has negative intersection against $g^{*}\alpha$. However, since $X$ is globally $F$-split, Lemma \ref{lemm:relcanoneff} shows that $K_{Z/X}$ has non-negative intersection against $g^{*}\alpha$, yielding a contradiction.  We conclude that $c_{1}(Q) \cdot \alpha \geq 0$ for every nef curve class $\alpha$.
\end{proof}

Theorem \ref{theo:weakpqp} indicates that the quotients with numerically trivial first Chern classes are ``boundary'' cases and thus are particularly interesting.  We next discuss the structure of such quotients.

\begin{corollary} \label{coro:maxtrivquotient}
    Suppose $X$ is a projective globally $F$-split variety.  There is a unique torsion-free quotient $\phi: T_{X} \to Q_{max}$ with the following properties:
    \begin{enumerate}
        \item $c_{1}(Q_{max})$ is a numerically trivial Weil divisor.
        \item If $\psi: T_{X} \to Q$ is any torsion-free quotient such that $c_{1}(Q)$ is numerically trivial, then $\psi$ factors through $\phi$.
    \end{enumerate}
\end{corollary}

\begin{proof}
    Let $\alpha$ be any curve class in the interior of $\Nef^{n-1}(X)$.  If $\mu^{min}(T_{X}) > 0$, then $Q_{max} = 0$.  If $\mu^{min}(T_{X}) = 0$, let $Q_{max}$ be the quotient of $T_{X}$ by the maximal positive subsheaf, i.e.~the next-to-last term $\mathcal{F}_{t-1}$ in the Harder-Narasimhan filtration.  Since $c_{1}(Q_{max})$ is pseudo-effective by Theorem \ref{theo:weakpqp}, the property $c_{1}(Q_{max}) \cdot \alpha = 0$ shows that $c_{1}(Q_{max})$ is numerically trivial.  Since $\mu^{min}(\mathcal{F}_{t-1}) > 0$, the subsheaf $\mathcal{F}_{t-1}$ must lie in the kernel of any other torsion-free quotient $\psi: T_{X} \to Q$ such that $c_{1}(Q)$ is numerically trivial.  This finishes the proof.
\end{proof}

In the situation of Corollary \ref{coro:maxtrivquotient}, we will call $Q_{max}$ the maximal trivial quotient of $T_{X}$.  Note that $Q_{max}$ is ``trivial'' in a strong sense -- it is semistable of slope $0$ with respect to every nef curve class.  However, the examples of \cite{Kollar95} show that $Q_{max}$ does not need to be numerically flat or essentially finite (see Example \ref{exam:kollar}).

A globally $F$-regular variety $X$ has the nicest geometry when $Q_{max}$ is zero.  This property is captured by the following definition.

\begin{definition} \label{defi:pqp}
    Let $X$ be a normal projective variety.  We say that $T_{X}$ has the positive quotient property (PQP) if for every nontrivial torsion-free quotient $Q$ of $T_{X}$ the first Chern class $c_1(Q)$ is pseudo-effective and not numerically trivial.
\end{definition}

This definition plays a key technical role in the paper.  For example, it is connected to the existence of rational curves and the non-existence of pluridifferential forms by the recent work of \cite{JLR25}.  In the next section, we will study how the PQP for the tangent bundle interacts with purely inseparable covers.

\section{Crepant covers} \label{sect:crepantcovers}

In Lemma \ref{lemm:relcaneffstrong} we saw that if $X$ is globally $F$-regular and $f: Y \to X$ is a finite cover then $K_{Y/X}$ is effective.  Recall from Definition \ref{defi:crepantcover} that a finite cover is crepant if $K_{Y/X} \sim 0$; these play a particularly important role in determining the geometry of $X$.

The goal of this section is to study the crepant covers $f: Y \to X$.  Since the separable crepant covers are controlled by the \'etale fundamental group of the smooth locus of $X$, our main focus is the purely inseparable covers.

\subsection{Purely inseparable crepant covers} \label{sect:allispicc}

Our first task is to prove a couple structural results for purely inseparable crepant covers.  We start by showing that crepancy is equivalent to a weaker hypothesis for purely inseparable covers of globally $F$-split varieties.

\begin{corollary} \label{coro:numtriv}
Let $X$ be a projective globally $F$-split variety.  Suppose that $f: Y \to X$ is a purely inseparable finite cover such that the divisor $K_{Y/X}$ is numerically trivial.  Then $K_{Y/X}$ is linearly equivalent to $0$.
\end{corollary} 

\begin{proof}
    Lemma \ref{lemm:relcanoneff} shows that $K_{Y/X}$ is linearly equivalent to an effective divisor.  The only effective divisor that is numerically equivalent to $0$ is the zero divisor.
\end{proof}

We next prove a basic fact about $F$-splitting and crepant covers.

\begin{proposition} \label{prop:foliationcoverfsplitting}
    Let $X$ be a projective globally $F$-split variety.  Suppose that $f: Y \to X$ is a purely inseparable crepant cover.  Then
    \begin{enumerate}
        \item $Y$ is globally $F$-split.
        \item If $X$ is globally $F$-regular, then so is $Y$.
    \end{enumerate}
\end{proposition}

\begin{proof}
    (1) Choose a positive integer $e$ such that $F^{e}_{rel}: X^{(-e)} \to X$ factors through $g: X^{(-e)} \to Y$.  Let $U \subset X$ and $V \subset Y$ be big smooth open subsets such that $f^{-1}U = V$.  By Grothendieck duality for the map $f$, we have
    \begin{equation*}
        \mathcal{H}om(f_{*}(g_{*}\mathcal{O}_{U^{(-e)}}),\mathcal{O}_{U}) = f_{*}\mathcal{H}om(g_{*}\mathcal{O}_{U^{(-e)}},\omega_{V/U}).
    \end{equation*}
    Since $X$ is $F$-split, we have a splitting on the left.  This corresponds to a homomorphism $\phi: g_{*}\mathcal{O}_{U^{(-e)}} \to \omega_{V/U} \cong \mathcal{O}_{V}$.  Since $H^{0}(V,\mathcal{O}_{V}) \cong \mathbb{K}$, $\phi(1)$ is a unit; after modifying $\phi$ by rescaling by this constant, we see that it is a splitting of the map $\mathcal{O}_{V} \to g_{*}\mathcal{O}_{U^{(-e)}}$. As the pushforward $g_{*}\mathcal{O}_{X^{(-e)}}$ is reflexive, by taking double duals we get a splitting of the map $\mathcal{O}_{Y} \to g_{*}\mathcal{O}_{X^{(-e)}}$.  Since $X^{(-e)}$ is globally $F$-split by Lemma \ref{lemm:frobeniustwistandsplitting}, so is $Y$ by \cite[Proposition 1.6]{HS11}.

    (2) This follows from the same argument as (1).
\end{proof}

We next show that purely inseparable crepant covers factor in a particularly simple way.  We will need the following lemma.

\begin{lemma} \label{lemm:negativityofpushox2}
    Let $X$ be a projective globally $F$-split variety and let $\alpha$ be a nef curve class on $X$.  Define $\mathcal{G}_{e} = F^{e}_{rel,*}\mathcal{O}_{X^{(-e)}}/F^{e-1}_{rel,*}\mathcal{O}_{X^{(1-e)}}$.  Then
    \begin{equation*}
        \mu^{max}_{F_{rel}^{*}\alpha}(F_{rel}^{\widehat{*}}\mathcal{G}_{e}) \leq 0 \qquad \textrm{and} \qquad \mu^{max}_{\alpha}(\mathcal{G}_{e}) \leq 0.
    \end{equation*}
\end{lemma}

\begin{proof}  
    It suffices to prove the first inequality as the second is an immediate consequence.  The proof is by induction on $e$. 
    
    We start with the base case $e=1$.  Let $U \subset X$ denote the big open locus where $X$ is smooth.  \cite[Theorem 3.7]{Sun08} constructs a filtration of $F_{rel}^{\widehat{*}}F_{rel,*}\mathcal{O}_{U^{(-1)}}$ whose $\ell$th graded piece is isomorphic to $T^{\ell}(\Omega_{U^{(-1)}})$ where $T^{\ell}$ denotes the $\ell$th truncated symmetric power functor (as defined in \cite[Definition~3.4]{Sun08}).  By definition, $T^{\ell}(\Omega_{U^{(-1)}})$ is a subsheaf of $\Omega^{\otimes \ell}_{U^{(-1)}}$.  Since $\mu_{F_{rel}^{*}\alpha}^{max}(\widehat{\Omega_{X^{(-1)}}}) \leq 0$ by Theorem \ref{theo:weakpqp}, we conclude from Corollary \ref{coro:fsplitstronglyss} that the same property holds for the double dual of the $\ell$th truncated symmetric power $T^{\widehat{\ell}}(\Omega_{X^{(-1)}})$, and hence also for $F_{rel}^{\widehat{*}}F_{rel,*}\mathcal{O}_{X^{(-1)}}$.  Finally, since $\mathcal{O}_{X^{(-1)}}$ has slope $0$, the slope inequality descends to  $F_{rel}^{\widehat{*}}F_{rel,*}\mathcal{O}_{X^{(-1)}}/ \mathcal{O}_{X^{(-1)}} \cong F_{rel}^{\widehat{*}}(F_{rel,*}\mathcal{O}_{X^{(-1)}}/\mathcal{O}_{X})$.

    Now we prove the induction step.  We let $\varphi: X^{(-1)} \to X$ denote the $\mathbb{K}$-semilinear isomorphism obtained by base change from the Frobenius map on $\mathbb{K}$.  Then we have the identity $\mathcal{G}_{e} = F_{rel,*}\varphi^{*}\mathcal{G}_{e-1}$.  Let $U_{e} \subset X$ denote the big open locus where $X$ is smooth and $\mathcal{G}_{i}$ is locally free for every $i \leq e$.  On the open set $U_{e}$, \cite[Theorem 3.7]{Sun08} constructs a filtration of $F_{rel}^{*}\mathcal{G}_{e}$ whose graded pieces have the form $T^{\ell}(\Omega_{U_{e}^{(-1)}}) \otimes \varphi^{*}\mathcal{G}_{e-1}|_{U_{e}}$.  By the same argument as before, $\mu_{F_{rel}^{*}\alpha}^{max}(T^{\widehat{\ell}}(\Omega_{X^{(-1)}})) \leq 0$.  We claim that $\mu^{max}_{F^{*}_{rel}\alpha}(\varphi^{*}\widehat{\mathcal{G}_{e-1}}) \leq 0$ as well.  Indeed, since $\varphi$ is $\mathbb{K}$-semilinear and $\mathbb{K}$ is algebraically closed, $\varphi$ identifies the degree function on $0$-cycles for $X$ and $X^{(-1)}$ and hence preserves intersection products.  Since $F^{*}_{rel}\alpha$ is proportional to the $\varphi^{-1}$-pullback of $\alpha$, the claim follows from the induction hypothesis.  The strong semistability property in Corollary \ref{coro:fsplitstronglyss} implies by \cite[Corollary 5.18]{JLR26} that maximum slopes behave well under reflexive tensor products.  Thus $\mu^{max}_{F^{*}_{rel}\alpha}(F_{rel}^{\widehat{*}}\mathcal{G}_{e}) \leq 0$ as desired.
\end{proof}

\begin{lemma} \label{lemm:pushstructuressslope0}
    Let $X$ be a globally $F$-split variety.  Suppose that $f: Y \to X$ is a purely inseparable crepant cover.  Then $f_{*}\mathcal{O}_{Y}$ is a reflexive self-dual sheaf.
    
    Furthermore, suppose that $X$ is projective.  Then $f_{*}\mathcal{O}_{Y}$ is strongly semistable of slope $0$ with respect to every nef curve class.  Moreover, for each $k \geq 1$ the intersection $f_{*}\mathcal{O}_{Y} \cap F^{k}_{rel,*}\mathcal{O}_{X^{(-k)}}$ is strongly semistable of slope $0$ as well.
\end{lemma}

\begin{proof}
    Let $U \subset X$ and $V \subset Y$ be smooth big open subsets such that $f^{-1}U = V$.  Then Grothendieck duality gives an identification
    \begin{equation*}
        \mathcal{H}om(f_{*}(f^{*}\omega_{U}),\omega_{U}) \cong f_{*}\mathcal{H}om(f^{*}\omega_{U},\omega_{V}).
    \end{equation*}
    Since $\omega_{V/U} = \mathcal{O}_{V}$, we deduce that $\mathcal{H}om(f_{*}\mathcal{O}_{V},\mathcal{O}_{U}) \cong \mathcal{H}om(f_*(f^* \omega_U), \omega_U) \cong f_{*}\mathcal{O}_{V}$.  Since reflexiveness is preserved by finite surjective pushforward, $f_{*}\mathcal{O}_{Y}$ is reflexive and we conclude that it is self-dual.

    Next suppose that $X$ is projective.  Since $f_{*}\mathcal{O}_{Y}$ is self-dual, for any nef curve class $\alpha$ we find that $\mu(f_{*}\mathcal{O}_{Y}) = 0$.  Choose a positive integer $e$ such that $F^{e}_{rel}$ factors through $f$ and consider the filtration
    \begin{equation*}
        \mathcal{O}_{X} \subset F_{rel,*}\mathcal{O}_{X^{(-1)}} \subset \ldots \subset F^{e}_{rel,*}\mathcal{O}_{X^{(-e)}}
    \end{equation*}
    By Lemma \ref{lemm:negativityofpushox2} each successive quotient $\mathcal{G}_{k}$ in this filtration satisfies $\mu^{max}(\mathcal{G}_{k}) \leq 0$.  We claim that for the induced filtration on $f_{*}\mathcal{O}_{Y}$, each successive quotient is either the zero sheaf or is strongly semistable of slope $0$.  Indeed, if for some $k$ the sheaf 
    $$\big(f_{*}\mathcal{O}_{Y} \cap F^{k}_{rel,*}\mathcal{O}_{X^{(-k)}}\big)/\big(f_{*}\mathcal{O}_{Y} \cap F^{k-1}_{rel,*}\mathcal{O}_{X^{(1-k)}}\big)$$
    had negative slope, then the slope bound from Lemma \ref{lemm:negativityofpushox2} would imply that $f_{*}\mathcal{O}_{Y}$ has negative slope, a contradiction.  %As this subquotient is also a subsheaf of 
    Thus, when this subquotient is non-zero, its slope must achieve the maximal possible value for a subsheaf of $F^{k}_{rel,*}\mathcal{O}_{X^{(-k)}}/F^{k-1}_{rel,*}\mathcal{O}_{X^{(1-k)}}$, which is $0$.  Since it cannot have a destabilizing subsheaf, it is strongly semistable by Corollary \ref{coro:fsplitstronglyss}.
\end{proof}

\begin{proposition} \label{prop:allarePICC}
    Let $X$ be a projective globally $F$-split variety.  Suppose that $f: Y \to X$ is a purely inseparable crepant cover.  Then $f$ factors as a sequence of height $1$ purely inseparable crepant covers
    \begin{equation*}
        Y = Z_{k} \xrightarrow{g_{k}} Z_{k-1} \xrightarrow{g_{k-1}} \ldots \xrightarrow{g_{2}} Z_{1} \xrightarrow{g_{1}} Z_{0} = X
    \end{equation*}
\end{proposition}

\begin{proof}
    Fix a class $\alpha$ in the interior of $\Nef^{\dim X - 1}(X)$.  The proof is by induction on the height of $f$.  When $f$ has height $1$ the statement is clear.

    We next prove the induction step.  
    Consider $\mathcal{E} := f_{*}\mathcal{O}_{Y} \cap F_{rel,*}\mathcal{O}_{X^{(-1)}}$.  This sheaf strictly contains $\mathcal{O}_{X}$ since $f$ factors through some height $1$ purely inseparable cover of $X$.  
    Lemma \ref{lemm:pushstructuressslope0} shows that $\mathcal{E}$ is strongly semistable of slope $0$.  
    We may view $\mathcal{E}$ as a subsheaf of $F_{abs,*}\OO_X$ through the identity $F_{rel,*}\OO_{X^{(-1)}} \cong F_{abs,*}\OO_X$.  Observe that the composition $F_{abs}^{*}F_{abs,*}\mathcal{O}_{X} \to F_{abs}^{*}F_{abs,*}\Omega_{X} \to \Omega_{X}$ is $\OO_X$-linear because $F_{abs,*}$ of the exterior derivative is $\mathcal{O}_{X}$-linear.  Restricting this map to $F_{abs}^*\mathcal{E} \subset  F_{abs}^*F_{abs,*}\OO_X$ and taking double duals produces an $\OO_X$-linear map $F_{abs}^{\widehat{*}}\mathcal{E} \to \widehat{\Omega_{X}}$.  In combination with Theorem \ref{theo:weakpqp},
    Lemma \ref{lemm:pushstructuressslope0} shows that the image sheaf $\mathcal{G}$ must be semistable of slope $0$.  By Ekedahl's correspondence, $f$ factors through the height $1$ purely inseparable cover $g: Z \to X$ defined by the kernel $\mathcal{F}$ of the map $T_{X} \to \mathcal{G}^{\vee}$.  Since $T_X/\mathcal{F}$ is a subsheaf of $\mathcal{G}^\vee$, $\mu(T_X/\mathcal{F}) \geq 0$, and $\mathcal{G}^\vee$ is semistable of slope $0$, we conclude that $T_X/\mathcal{F}$ is also semistable of slope $0$.  By Corollary \ref{coro:numtriv} $g$ is a crepant cover.  As the induced map $Y \to Z$ is a purely inseparable crepant cover of a globally $F$-split variety with strictly smaller height, we conclude the desired statement by applying the induction assumption to $Y \to Z$.
\end{proof}

\subsection{Superior covers of globally $F$-regular varieties}

We will now prove the existence and uniqueness of superior covers, as defined in Definition \ref{def: superior cover}.

\begin{proposition} \label{prop:constructingpqpcover}
    Suppose $X$ is a projective globally $F$-regular variety.  There is a purely inseparable crepant cover $f: Z \to X$ such that $T_{Z}$ has the PQP.
\end{proposition}

\begin{proof}
   We construct a chain of height $1$ purely inseparable covers as follows.  Fix a nef curve class $\alpha$ in the interior of $\Nef^{\dim X-1}(X)$.  Since $X$ is globally $F$-regular, Corollary \ref{coro:pclosedcriterion} shows that its maximal positive subsheaf is a $p$-closed foliation.  By Lemma \ref{lemm:relCanonicalHeight1} and Corollary \ref{coro:numtriv} it defines a crepant height $1$ purely inseparable cover $g_{1}: Z_{1} \to X$.  By Proposition \ref{prop:foliationcoverfsplitting}.(2), $Z_{1}$ is also globally $F$-regular.  Thus we may continue this procedure inductively to obtain a chain
   \begin{equation*}
       \ldots \xrightarrow{g_{3}} Z_{2} \xrightarrow{g_{2}} Z_{1} \xrightarrow{g_{1}} X.
   \end{equation*}
   Suppose we fix an ample Cartier divisor $H$ on $X$.  Denote the pullback of $H$ to $Z_{i}$ by $H_{i}$ and let $S_{i}$ denote the section ring $\oplus_{m \geq 0} H^{0}(Z_{i},mH_{i})$.  
   The affine cone $C_i = \Spec S_i$ is globally $F$-regular by \cite[Proposition~5.3]{SS10}.  Moreover, on the complement of the origin the natural map $f_i : C_i \to C_{i-1}$ is a flat base change of $g_i$.  Hence, this map is also a purely inseparable crepant cover, and from Lemma \ref{lemm:pushstructuressslope0} we obtain that $S_i \cong \Hom_{S_{i-1}}(S_i, S_{i-1})$ as $S_i$-modules.

   Consider a splitting $g_{i*}\OO_{Z_i} \to \OO_{Z_{i-1}}$.  By the projection formula, this induces a graded $S_{i-1}$-module homomorphism $T_i : S_i \to S_{i-1}$ with $T_i(1) = 1$.  In particular, $T_{i}$ is a generator of the rank $1$ free $S_{i}$-module $\Hom_{S_{i-1}}(S_i, S_{i-1})$.  
   
   Localize $S_i$ and $S_{i-1}$ at the irrelevant ideals to obtain local rings $(\widetilde{S}_{i},\mathfrak{m}_{i})$ and $(\widetilde{S}_{i-1},\mathfrak{m}_{i-1})$.  As $S_{i-1} \to S_{i}$ is purely inseparable, the irrelevant ideal of $S_{i}$ is the only prime ideal lying over the irrelevant ideal of $S_{i-1}$ so that the induced map $\widetilde{S}_{i-1} \to \widetilde{S}_{i}$ is a finite local homomorphism.  Since $S_{i}$ is a finitely generated $S_{i-1}$-module, $\Hom$ commutes with localization, so the localized map $\widetilde{T}_{i}$ is still a generator of $\widetilde{S}_{i}$-module $\Hom_{\widetilde{S}_{i-1}}(\widetilde{S}_{i},\widetilde{S}_{i-1})$.  Since $\widetilde{T}_{i}$ satisfies $\widetilde{T}_{i}(1) = 1$, it is surjective.  Finally since $T_{i}$ preserves the grading we have $\widetilde{T}_{i}(\mathfrak{m}_{i}) \subset \mathfrak{m}_{i-1}$.  We conclude that we are in the setting of \cite[Theorem~4.8]{CR22}.  The theorem shows that the $F$-signature of $\widetilde{S}_{i}$ is at least $p$ times the $F$-signature of $\widetilde{S}_{i-1}$ unless $g_{i}$ is an isomorphism.  Since the $F$-signature of a strongly $F$-regular ring is positive and is bounded above by $1$, we conclude that our inductive procedure must stop with the desired cover $Z$.
\end{proof}

\begin{remark}
    \cite{Shen10} shows that a projective variety that is freely rationally connected will admit a purely inseparable cover that is separably rationally connected.  This result is similar in spirit to Proposition \ref{prop:constructingpqpcover} and we were motivated by Shen's work.
\end{remark}

\begin{proposition} \label{prop:allfactorthrough}
    Let $X$ be a projective globally $F$-regular variety.  Suppose that $f: Y \to X$  is a purely inseparable crepant cover such that $T_{Y}$ has the PQP.  Then every other purely inseparable crepant cover $g: Z \to X$ has the property that $f$ factors through $g$.
\end{proposition}

\begin{proof}
    Let $W$ denote the normalization of $(Y \times_{X} Z)_{red}$ and consider the finite morphism $h: W \to Y$.  We claim that $h$ is again a purely inseparable crepant cover.  Indeed, suppose first that both $f$ and $g$ have height $1$.  Then $f$ and $g$ are the covers associated to $p$-closed foliations $\mathcal{F}_{Y},\mathcal{F}_{Z} \subset T_{X}$.  Thus $W \to X$ is the height $1$ cover induced by $\mathcal{F}_{Y} \cap \mathcal{F}_{Z}$.  Since $T_{X}/(\mathcal{F}_{Y} \cap \mathcal{F}_{Z})$ has numerically trivial first Chern class by \cite[Lemma 4.27]{JLR26}, $W \to X$ is crepant by Corollary \ref{coro:numtriv}.  Thus $h: W \to Y$ is also crepant.  In the general case, we apply Proposition \ref{prop:allarePICC} to $f$ and $g$ to write them as compositions of height $1$ purely inseparable covers.  We then apply the previous argument inductively to realize $h$ as a sequence of height $1$ purely inseparable crepant covers.

    Since $T_{Y}$ has the positive quotient property, by Proposition \ref{prop:allarePICC} $h$ must be an isomorphism.  This implies that as subfields of $K(X)^{1/p^{\infty}}$ we have $K(Z) \cdot K(Y) = K(Y)$, i.e.~$K(Z) \subset K(Y)$.  The claim follows.
\end{proof}

Putting everything together, the following theorem identifies a canonical purely inseparable crepant cover $f: Y \to X$ of a projective globally $F$-regular variety $X$.  This map $f$ will of course be the superior cover; we will verify this in Section \ref{sect:superior} once we connect the PQP of the tangent bundle with the vanishing of pluridifferential forms.

\begin{theorem} \label{theo:pqpcanonicalcover}
    Let $X$ be a projective globally $F$-regular variety.  There is a unique purely inseparable crepant cover $f: Y \to X$ such that $T_{Y}$ has the PQP.  Furthermore:
    \begin{itemize}
        \item $Y$ is globally $F$-regular.
        \item $f$ factors through every purely inseparable crepant cover of $X$.
    \end{itemize}
\end{theorem}

\begin{proof}
    The existence of a purely inseparable crepant cover $f: Y \to X$ whose tangent bundle satisfies the PQP is established by Proposition \ref{prop:constructingpqpcover}.  The variety $Y$ is globally $F$-regular by Proposition \ref{prop:foliationcoverfsplitting}.  The desired factoring property is established in Proposition \ref{prop:allfactorthrough}; this also implies the uniqueness of $f$.
\end{proof}

\subsection{Finiteness of purely inseparable crepant covers}

We next turn to a proof of Theorem \ref{theo:finiteness}, which will use Tannakian categories.  Similar ideas have been used to study essentially finite sheaves (via Nori's fundamental group scheme; see \cite{Nori82}) and numerically flat sheaves (via Langer's S-fundamental group scheme; see \cite{Langer11}).  For a given globally $F$-regular variety $X$, $Q_{max}$ might be neither essentially finite nor numerically flat.  Nevertheless, the weaker form of triviality satisfied by $Q_{max}$ still allows us to leverage the additional structure provided by Tannakian categories.

Suppose $X$ is a projective globally $F$-split variety and $\alpha$ is a curve class in the interior of $\Nef^{\dim X - 1}(X)$.  Let $\mathcal{C}$ denote the category of reflexive sheaves that are semistable of slope $0$.  We equip $\mathcal{C}$ with the following structures:
\begin{itemize}
    \item $\mathcal{C}$ is an abelian category, if we include double duals when defining cokernels.
    \item $\mathcal{C}$ carries a tensor product structure $\widehat{\otimes}$.  (Note that this operation preserves reflexiveness and preserves semistability by Corollary \ref{coro:fsplitstronglyss}.)
\end{itemize}
It is easy to check that these give $\mathcal{C}$ the structure of a Tannakian category, where the fiber functor is obtained by taking stalks at the generic point of $X$.

The following result was communicated to us by ChatGPT.  It shows that the number of height $1$ purely inseparable crepant covers is finite; the key idea is to linearize the problem using the Tannakian structure.  The overall structure of the proof is very similar to \cite[Corollary 3.5]{AE11} and \cite[Proposition 13]{BS11} which address essentially finite bundles.  The Grassmannian argument is based on \cite[Corollary 5.(iii)]{BS11} and ChatGPT identifies \cite{GP26} as an important precursor for involving the algebra structure.

\begin{theorem} \label{theo:finmanyheight1}
    Let $X$ be a projective globally $F$-split variety.  There are only finitely many height $1$ purely inseparable crepant covers $f: Y \to X$ up to isomorphism.
\end{theorem}

\begin{proof}
    Throughout we fix a class $\alpha$ in the interior of the nef cone of curves and a splitting $\phi: F_{abs,*}\mathcal{O}_{X} \to \mathcal{O}_{X}$.
    
    We let $\mathfrak{S}$ denote the set of unital $\mathcal{O}_{X}$-subalgebras $\mathcal{B} \subset F_{abs,*}\mathcal{O}_{X}$ which are reflexive and semistable of slope $0$.
    Note that if $f: Y \to X$ is a height $1$ purely inseparable crepant cover, then $f_{*}\mathcal{O}_{Y} \in \mathfrak{S}$ by Lemma \ref{lemm:pushstructuressslope0}.  Since $f$ is determined by $f_{*}\mathcal{O}_{Y}$, it suffices to prove that $\mathfrak{S}$ is a finite set.

    The map $\phi$ induces a canonical identification of each algebra $\mathcal{B} \in \mathfrak{S}$ with its dual.  Indeed, consider the map $\rho_{\mathcal{B}}: \mathcal{B} \to \mathcal{B}^{\vee}$ sending $b \mapsto \phi(b \cdot -)$.  We claim that $\rho_{\mathcal{B}}$ is an isomorphism.  First note that it is an isomorphism at the generic point: for any non-zero $b$, we have $\phi(b \cdot b^{-1}) = 1$ so that the pairing is non-degenerate.  Since both $\mathcal{B}$ and $\mathcal{B}^{\vee}$ are semistable of slope $0$ and since $\alpha$ lies in the interior of the nef cone, any map that is generically an isomorphism must also be an isomorphism in codimension $1$.  Finally, since both sheaves are reflexive this isomorphism must extend over all of $X$.
    
    By Lemma \ref{lemm:negativityofpushox2} we have $\mu^{max}(F_{abs,*}\mathcal{O}_{X}) = 0$.  Let $\mathcal{A}$ be the maximal slope $0$ subsheaf of $F_{abs,*}\mathcal{O}_{X}$.  It is clear that $\mathcal{O}_{X} \subset \mathcal{A}$.  Furthermore since $\mathcal{A} \widehat{\otimes} \mathcal{A}$ is semistable of slope $0$ by Corollary \ref{coro:fsplitstronglyss} and \cite[Corollary 5.18]{JLR26}, its image in $F_{abs,*}\mathcal{O}_{X}$ under the multiplication map must lie in $\mathcal{A}$ so that $\mathcal{A}$ is a unital algebra.  Note also that $\mathcal{A}$ is a saturated subsheaf of the reflexive sheaf $F_{abs,*}\mathcal{O}_{X}$, hence reflexive.  We conclude that $\mathcal{A} \in \mathfrak{S}$.
    
    We define $\mathcal{C}_{\mathcal{A}} \subset \mathcal{C}$ to be the tensor subcategory generated by subquotients of $\mathcal{A}$.  By \cite[Corollaire 6.20]{Deligne90}, $\mathcal{C}_{\mathcal{A}}$ has a fiber functor defined over a finite extension of the ground field $\mathbb{K}$.  Since $\mathbb{K}$ is algebraically closed, we obtain a fiber functor over $\mathbb{K}$, yielding an equivalence $\psi: \mathcal{C}_{\mathcal{A}} \xrightarrow{\cong} \mathrm{Rep}_{\mathbb{K}}(G)$ for some affine group scheme $G$ of finite type.  We let $V$ denote $\psi(\mathcal{A})$.  The algebra structure on $\mathcal{A}$ induces an algebra structure $V \cdot V \to V$.  Furthermore $V$ carries a bilinear form $q$ associated with the isomorphism $\rho_{\mathcal{A}}$.

    Every element of $\mathfrak{S}$ is contained in $\mathcal{C}_{\mathcal{A}}$.  We can describe their images under $\psi$ as follows.  Consider the Grassmannian of rank $r$ subspaces $\mathrm{Gr}(r,V)$.  Define $Z_{r} \subset \mathrm{Gr}(r,V)$ to be the sublocus
    \begin{equation*}
    Z_{r} = \{ W \subset V \, |  \, 1_{V} \in W, \, W \cdot W \subset W, \, W \textrm{ is }G\textrm{-stable} \}.
    \end{equation*}
    As we vary over all $r$ the elements of $Z_{r}$ correspond exactly to the elements of $\mathfrak{S}$.  Note that $Z_{r}$ is defined by closed conditions and thus defines a projective subscheme of $\mathrm{Gr}(r,V)$.
    
    Over $\mathrm{Gr}(r,V)$, $q$ induces a map of tautological bundles $q_{\mathcal{S}}: \mathcal{S} \to \mathcal{S}^{\vee}$.  Its determinant is a non-zero section $\sigma \in H^{0}(\mathrm{Gr}(r,V),\mathcal{O}(2))$.  
    We verified earlier that the restriction of $q$ to $\psi(\mathcal{B})$ is a non-degenerate bilinear form for every $\mathcal{B} \in \mathfrak{S}$.  Thus $\sigma|_{Z_{r}}$ is nowhere vanishing.  Since $\mathcal{O}(2)$ is ample, we must have $\dim Z_{r} = 0$.  By letting $r$ vary, we obtain the desired finiteness statement. 
\end{proof}

\begin{example}
    An $F$-splitting assumption is necessary in Theorem \ref{theo:finmanyheight1}.  For an abelian variety $X$, the height $1$ purely inseparable crepant covers correspond to height $1$ subgroup schemes of $X$.  There are only finitely many such covers when $X$ is ordinary, but there can be infinitely many when $X$ is supersingular.
\end{example}

\begin{proof}[Proof of Theorem \ref{theo:finiteness}:]
    Consider the directed graph $T$ whose nodes are purely inseparable crepant covers $f: Y \to X$ (up to isomorphism over $X$) and whose directed edges represent height $1$ purely inseparable maps $g: Y_{1} \to Y_{2}$ yielding a commutative diagram
    \[
    \xymatrix@R=0.3pc{
        Y_{1} \ar[dd]^{g} \ar[dr] &  \\
        & X \\
        Y_{2} \ar[ur] &
    }
    \]
    Proposition \ref{prop:allarePICC} shows that $T$ is connected.  Theorem \ref{theo:finmanyheight1} shows that every node receives only finitely many edges.  Theorem \ref{theo:pqpcanonicalcover} shows that the graph has finite diameter.  Since this graph has a unique terminal element, the desired finiteness follows.
\end{proof}

\subsection{Separable crepant covers}

As mentioned in the introduction, our work also allows us to give a new proof of a result of \cite{CRST18}.

\begin{theorem}[{\cite[Corollary 5.7]{CRST18}}] \label{theo:finiteetalefund}
Let $X$ be a projective globally $F$-regular variety.  Then the \'etale fundamental group of the smooth locus of $X$ is finite.
\end{theorem}

\begin{proof}
    Let $g: Z \to X$ be the purely inseparable finite cover constructed in Proposition \ref{prop:constructingpqpcover}.  \cite[Theorems 1.7 and 1.8]{JLR25} show that the \'etale fundamental group of the smooth locus of $Z$ is finite.  By purity we deduce finiteness of the \'etale fundamental group for any big open subset of the smooth locus of $Z$, and thus also for any big open subset $V \subset Z$.

    Let $U \subset X$ be the smooth locus and set $V = g^{-1}U$.  The map $g: V \to U$ is a universal homeomorphism.  By \cite[Expos\'e IX, Th\'eor\`eme 4.10]{SGA1}, we have an isomorphism
    \begin{equation*}
        g_{*}: \pi_{1}^{et}(V) \xrightarrow{\cong} \pi_{1}^{et}(U).
    \end{equation*}
    We conclude finiteness of $\pi_{1}^{et}(U)$ from the previous paragraph.
\end{proof}

\cite[Lemma 11.1]{PZ26} shows that a quasi-\'etale cover of a globally $F$-regular variety is again globally $F$-regular.  Thus our results can be applied inductively to analyze any crepant cover $f: Y \to X$ that can be written as a composition of separable and purely inseparable crepant covers. However, we do not know the answer to the following natural question:

\begin{question}
    Let $X$ be a projective globally $F$-regular variety.  Can every crepant cover of $X$ be written as a composition of purely inseparable crepant covers and separable crepant covers?
\end{question}

\section{Superior globally $F$-regular varieties} \label{sect:superior}

The following result connects the PQP with the other desirable geometric properties discussed in the introduction.  In particular, it shows that a projective globally $F$-regular variety $X$ will be superior if and only if $T_{X}$ has the PQP.

\begin{theorem} \label{theo:bodyequiv}
    Let $X$ be a projective globally $F$-regular variety.  Then the following are equivalent:
    \begin{enumerate}
        \item $H^{0}(X,\widehat{\Omega_{X}^{\otimes j}})=0$ for every $j \geq 1$.
        \item $T_{X}$ has the positive quotient property.
        \item $X$ has no non-trivial purely inseparable crepant covers.
        \item $X$ admits a very free curve in its smooth locus (possibly of positive genus).
    \end{enumerate}
\end{theorem}

\begin{proof}
    (1) $\Rightarrow$ (2): Let $T_{X} \to Q_{max}$ be the quotient of Corollary \ref{coro:maxtrivquotient} and suppose for a contradiction that the rank $r$ of $Q_{max}$ is greater than $0$.  Let $f: Y \to X$ denote the crepant cover corresponding to $Q_{max}$.  Lemma \ref{lemm:relCanonicalHeight1} and Corollary \ref{coro:numtriv} show that $(p-1) \cdot f^{*}c_{1}(Q_{max}) \sim 0$.  Pushing forward, we see $(p-1)\deg(f)c_{1}(Q_{max}) \sim 0$.  Set $m = (p-1)\deg(f)$ and let $U \subset X$ be the smooth big open subset where $Q_{max}$ is locally free.  Then we get a surjection
    \begin{equation*}
    T_{U}^{\otimes mr} \to Q_{max}|_{U}^{\otimes mr} \to \left(\bigwedge\nolimits^{r} Q_{max}|_{U} \right)^{\otimes m}.
    \end{equation*}
    Taking duals of the reflexive closures on $X$, we get a section of $\widehat{\Omega_{X}^{\otimes mr}}$.

    (2) $\Leftrightarrow$ (3): Follows from Theorem \ref{theo:pqpcanonicalcover}. 

    (2) $\Rightarrow$ (4):  This follows from \cite[Theorem 1.7]{JLR25}.

    (4) $\Rightarrow$ (1): Immediate.
\end{proof}

\begin{proof}[Proof of Theorem \ref{theo:canonicalcover}:]
Follows from Theorem \ref{theo:bodyequiv} and Theorem \ref{theo:pqpcanonicalcover}.
\end{proof}

\begin{proof}[Proof of Theorem \ref{theo:introequiv}:]
Follows from Theorem \ref{theo:bodyequiv}.
\end{proof}

\begin{proof}[Proof of Theorem \ref{theo:maintheorem2}:]
Follows from Theorem \ref{theo:bodyequiv} and Theorem \ref{theo:pqpcanonicalcover}.
\end{proof}

\subsection{Smooth globally $F$-regular varieties} \label{sect:smoothones}

Smooth projective globally $F$-regular varieties share some of the key properties of smooth Fano varieties in characteristic $0$.  For example, they are ``simply connected'' in a strong sense: they have a trivial Nori fundamental group scheme (\cite[Theorem (A)]{KV25}) and a trivial $S$-fundamental group scheme (\cite[Theorem 1.2]{EM11}).  They also satisfy vanishing of certain differential forms: \cite[Theorem 10.2]{KV25} shows that $H^{0}(X,\Omega^{1}_{X}) = 0$.

Under a smoothness assumption, we can prove a stronger version of Theorem \ref{theo:bodyequiv}.

\begin{proposition} \label{prop:smoothcase}
    Let $X$ be a smooth projective globally $F$-regular variety.  Then the following conditions are equivalent:
    \begin{enumerate}
        \item $H^{0}(X,\wedge^{j} \Omega_{X})=0$ for every $j > 0$.
        \item $T_{X}$ has the positive quotient property.
        \item $X$ has no non-trivial separable or purely inseparable crepant covers.
        \item $X$ admits a very free curve (possibly of positive genus).
    \end{enumerate}
\end{proposition}

\begin{proof}
    Most of the steps are the same as in the proof of Theorem \ref{theo:bodyequiv}.  We will only address the steps that are different.

    (1) $\Rightarrow$ (2): Let $T_{X} \to Q_{max}$ be the quotient of Corollary \ref{coro:maxtrivquotient} and suppose for a contradiction that the rank $r$ of $Q_{max}$ is greater than $0$.  \cite[Corollary 5.4]{CR22} shows that $\Pic(X)$ has no torsion so that $c_{1}(Q_{max}) \sim 0$.  Then the map $\widehat{\wedge}^{r} Q_{max}^{\vee} \to \wedge^{r} \Omega^{1}_{X}$ constructs a section, contradicting our vanishing assumption.

    (2) $\Rightarrow$ (3): This requires the additional well-known fact that the \'etale fundamental group of $X$ is trivial; see e.g.~\cite[Theorem (A)]{KV25}.
\end{proof}

Based on the evidence mentioned above, it is natural to wonder whether smooth projective globally $F$-regular varieties are always superior.  We conjecture that this is the case.

\begin{conjecture} \label{conj:smoothcase}
    Let $X$ be a smooth projective globally $F$-regular variety.  Then $T_{X}$ has the positive quotient property.
\end{conjecture}

\section{Rational curves} \label{sect:rationalcurves}

We next turn to the study of rational curves on projective globally $F$-regular varieties.  Prior work has focused on low dimensions and on special classes of varieties.  For example, rational chain connectedness has been established for several general classes of projective threefolds, including globally $F$-regular pairs with $p \geq 11$ (\cite[Theorem 0.1]{GLPSTZ15}) and klt Fano type varieties with $p>5$ (\cite[Theorem 1.5]{GNT19}).  Recently \cite{OPTWX26} has extended these results to mixed characteristic for certain globally $+$-regular pairs of (relative) dimension $3$.   Additionally, there has been substantial work on separable rational connectedness for hypersurfaces (\cite{BV17}, \cite{ST19}, \cite{BS23}, \cite{Zhu24}, \cite{BMR}), complete intersections (\cite{CZ14}, \cite{Tian15}, \cite{CR19}, \cite{STZ22}, \cite{CS22}, \cite{BVY24}), and varieties with nef tangent bundle (\cite{KW23}). 

We first recall several different ways of identifying when a variety carries ``many'' rational curves.

\begin{definition}
    Let $X$ be a quasiprojective variety.  We say that:
    \begin{itemize}
        \item $X$ is rationally chain connected (RCC), if any two general points of $X$ are connected by a chain of proper rational curves in $X$;
        \item $X$ is rationally connected (RC), if there is a family of proper rational curves $U \to M$ and a morphism $f: U \to X$ such that the induced map $U \times_{M} U \to X \times X$ is dominant;
        \item $X$ is separably rationally connected (SRC), if there is a smooth family of proper rational curves $U \to M$ and a morphism $f: U \to X$ such that the induced map $U \times_{M} U \to X \times X$ is dominant and generically smooth;
        \item $X$ is freely rationally connected (FRC), if there is a smooth family of proper rational curves $U \to M$ and a morphism $f: U \to X$ such that the general curve parametrized by $M$ is a free rational curve mapping to the smooth locus $X_{sm}$ and the induced map $U \times_{M} U \to X \times X$ is dominant;
        \item $X$ is smoothly separably rationally connected (SSRC), if the smooth locus $X_{sm}$ is SRC, or equivalently, if there is a very free rational curve in $X_{sm}$.
    \end{itemize}
\end{definition}
There are implications
\begin{equation*}
\begin{tikzcd}[column sep=1.5em, row sep=-0.5em]
    & X \textrm{ SRC} \arrow[dr, Rightarrow] & & \\
    X \textrm{ SSRC} \arrow[ur, Rightarrow] \arrow[dr, Rightarrow]
    & & X \textrm{ RC} \arrow[r, Rightarrow] & X \textrm{ RCC} \\
    & X \textrm{ FRC} \arrow[ur, Rightarrow] & &
\end{tikzcd}
\end{equation*}
but none of these conditions are equivalent in positive characteristic.

\subsection{Rational chain connectedness}

We are now ready to prove that projective globally $F$-regular varieties are RCC.

\begin{proof}[Proof of Theorem \ref{theo:maintheorem1}:]
    Let $f: Y \to X$ be the purely inseparable finite cover constructed in Theorem \ref{theo:pqpcanonicalcover}.  By \cite[Theorem 1.7]{JLR25} there is a smooth projective curve $C$ admitting a very free map $s: C \to Y$.

    Let $M \subset \Mor(C,Y)$ denote the irreducible component containing $s$ (which is unique since $s$ is a smooth point).  By \cite[Lemma 3.6]{LRT24} the two-fold evaluation map $C \times C \times M \to (Y \times C)^{\times 2}$ is dominant. 
    
    We claim that $Y$ is rationally chain connected.  If $g(C) = 0$ we are done.  Otherwise, let $h: C \to Y$ be a general map in $M$ sending two general points $p_{1},p_{2}$ in $C$ to two general points $y_{1},y_{2}$ in $Y$.  Since $H^{1}(C,h^{*}T_{Y}(-p_{1}-p_{2}))=0$, the natural map $H^{0}(C,N_{h}(-p_{1}-p_{2})) \to H^{1}(C,T_{C}(-p_{1}-p_{2}))$ is surjective.  This shows that in the Kontsevich space of pointed stable maps $\overline{\mathcal{M}}_{g(C),2}(Y)$, the irreducible component of $ev_{1}^{-1}(y_{1}) \cap ev_{2}^{-1}(y_{2})$ containing $h$ dominates $\overline{\mathcal{M}}_{g(C),2}$ under the forgetful map.  In particular, $h$ degenerates in the Kontsevich space to a morphism $h': C' \to Y$ whose domain is a nodal curve with geometrically rational components such that $h'(C')$ still contains $y_{1}$ and $y_{2}$. This proves that $Y$, and hence also $X$, is rationally chain connected.
\end{proof}

\begin{remark}
    The proof shows that every projective globally $F$-regular variety is connected by genus $g(C)$ curves in the sense of \cite{Gounelas16}: through any two general points (or indeed, if we use $r$-free curves in place of very free curves, any general finite set) there is an irreducible curve of genus $g(C)$.  
\end{remark}

\begin{example}
    \cite{KW23} shows that if $X$ is a smooth projective rationally chain connected variety such that $T_{X}$ is nef then $X$ admits a very free rational curve.  Our work shows that the rationally chain connected assumption can be replaced by  assuming either (1) $X$ is globally $F$-regular, or (2) $T_{X}$ has the PQP.
\end{example}

\subsection{Stronger conjectures}

It is natural to wonder whether one can strengthen the relationship between rational curves and globally $F$-regular varieties presented in Theorem \ref{theo:maintheorem1}.  \cite{KM99} conjectures that a klt Fano variety in characteristic $0$ will contain a very free rational curve in its smooth locus.  This property does not hold for every projective globally $F$-regular variety -- it is obstructed by the existence of pluridifferential forms.

We expect that for projective globally $F$-regular varieties being superior is equivalent to being SSRC. Recall that Theorem \ref{theo:introequiv} gives a list of equivalent conditions that characterize the superior globally $F$-regular varieties.  Each of these conditions is necessary for $X$ to be SSRC.  It is quite natural to expect them to be sufficient as well.

\begin{conjecture} \label{conj:pqpthenfree}
    Let $X$ be a superior globally $F$-regular variety.  Then there are very free rational curves in the smooth locus of $X$.
\end{conjecture}

\begin{example} \label{exam:introtoric}
    \cite{Smith00} shows that every projective toric variety is globally $F$-regular.  However, there are examples of toric varieties which do not carry any very free rational curves in their smooth locus.  We will characterize all such toric varieties in Proposition \ref{prop:toricexam}, verifying Conjecture \ref{conj:pqpthenfree} in the toric case. 
\end{example}

Via Theorem \ref{theo:canonicalcover}, Conjecture \ref{conj:pqpthenfree} would imply that a projective globally $F$-regular variety is rationally connected.  Based on Theorem \ref{theo:weakpqp} and the similarity with the results of \cite{Shen10}, we conjecture a stronger statement:

\begin{conjecture} \label{conj:freelyratconn}
    Let $X$ be a projective globally $F$-regular variety.  Then $X$ is freely rationally connected.
\end{conjecture}

Conjecture \ref{conj:smoothcase} predicts that every smooth projective globally $F$-regular variety is superior.  Combined with Conjecture \ref{conj:pqpthenfree}, we expect that Koll\'ar's question concerning very free curves on smooth Fano varieties has an affirmative answer under an additional globally $F$-regular assumption.

\begin{conjecture}[{\cite[Conjecture 3.7]{LT26}}]
    Let $X$ be a smooth projective globally $F$-regular variety.  Then $X$ admits a very free rational curve.
\end{conjecture}

\section{Examples}  \label{sect:examples}

Finally, we illustrate the general theory in a couple of examples.  In the first subsection we analyze toric varieties: we characterize the superior toric varieties and explain how to construct the superior cover using toric techniques.  We also verify that the superior condition is equivalent to being SSRC, verifying Conjecture \ref{conj:pqpthenfree} in the toric case.

In the second subsection we study the examples of \cite{Kollar95} where $Q_{max}$ has interesting properties.

\subsection{Toric varieties}

We next systematically analyze normal projective toric varieties.  Every such variety is globally $F$-regular by \cite[Proposition 6.4]{Smith00}.  Toric varieties are also SRC since they are birationally equivalent to projective space.  (In fact, using the torus action one can show that toric varieties are FRC.) The following proposition characterizes the SSRC condition.  (See also \cite{FM25,CLT26}.)

\begin{proposition} \label{prop:toricexam} 
    Suppose that $X$ is a normal projective toric variety defined by a fan $\Sigma \subset N$.  The following conditions are equivalent:
    \begin{enumerate}
        \item $X$ is SSRC. 
        \item The only torsion-free quotient of $T_{X}$ with numerically trivial first Chern class is the zero sheaf.
        \item The primitive generators $v_{\rho}$ of the rays $\rho \in \Sigma$ span $N \otimes \mathbb{F}_{p}$.
    \end{enumerate}
\end{proposition}

Note that condition (3) (and thus all the conditions) will hold whenever $X$ is smooth.

\begin{proof}
    (1) $\implies$ (2).  Clear.
    
    (2) $\implies$ (3).  Since the log tangent sheaf of $X$ is globally generated, every torsion-free quotient of $T_{X}$ will admit a globally generated subsheaf of full rank.  Thus, if a torsion-free quotient of $T_{X}$ has numerically trivial first Chern class then its reflexive hull must be a direct sum of $\mathcal{O}_{X}$.  In particular, $T_{X}$ will have such a quotient if and only if $\widehat{\Omega^{1}_{X}}$ admits a non-zero section.  In the standard exact sequence
    \begin{equation} \label{eq:torictangent}
        0 \to \widehat{\Omega^{1}_{X}} \to N^{\vee} \otimes \mathcal{O}_{X} \to \oplus_{\rho} \mathcal{O}_{D_{\rho}}
    \end{equation}
    the rightmost map sends $m \otimes f \mapsto (\langle m,v_{\rho} \rangle f|_{D_{\rho}})_{\rho}$.  Thus $\widehat{\Omega^{1}_{X}}$ will have a non-zero section precisely when the set of primitive generators $\{ v_{\rho} \}$ in $N$ fail to span $N \otimes \mathbb{F}_{p}$.

    (3) $\implies$ (1). Suppose that the generators $\{ v_{\rho} \}$ span $N \otimes \mathbb{F}_{p}$.  One can explicitly construct morphisms $s: \mathbb{P}^{1} \to X$ in the following way.  By projectivity, there is a linear combination $\sum_{\rho} a_{\rho}v_{\rho} = 0$ where each $a_{\rho}$ is a positive integer.  Let $T \subset X$ be the torus and for each ray $\rho$ of $\Sigma$, let $\varphi_{\rho}: \mathbb{G}_{m} \to T$ denote the corresponding one-parameter subgroup.  One can define a morphism $s: \mathbb{P}^{1} \to X$ by choosing different elements $\lambda_{\rho} \in \mathbb{K}$ and taking the closure of the rational map $s: \mathbb{G}_{m} \dashrightarrow T$ defined as
    \begin{equation*}
        s(t) = \prod_{\rho} \varphi_{\rho}(t-\lambda_{\rho})^{a_{\rho}}.
    \end{equation*}
    The image of $s$ will only meet the toric boundary at general points of codimension $1$ loci and thus will lie in the smooth locus of $X$.  Note that $s^{*}D_{\rho}$ is concentrated at $\lambda_{\rho}$ and has multiplicity $a_{\rho}$.  Thus by restricting Equation \eqref{eq:torictangent} we obtain
    \begin{equation*}
        0 \to s^{*}\widehat{\Omega^{1}_{X}} \to N^{\vee} \otimes \mathcal{O}_{\mathbb{P}^{1}} \to \oplus_{\rho} \mathcal{T}_{a_{\rho}} \to 0
    \end{equation*}
    where $\mathcal{T}_{a_{\rho}}$ is the structure sheaf of the length $a_{\rho}$-subscheme supported at $\lambda_{\rho}$.  Since the $v_{\rho}$ span $N \otimes \mathbb{F}_{p}$ and the rightmost map is defined by pairing against the various $v_{\rho}$, we see that $s$ is a very free rational curve, showing that $X$ is SSRC.
\end{proof}

Since the tangent bundle of the toric variety need not have the PQP, it is interesting to ask how to interpret the superior cover in the toric setting.  This is accomplished in the next construction.

\begin{construction} \label{cons:toricsuperior}
    Suppose that $X$ is a normal projective toric variety defined by a fan $\Sigma \subset N$.  The superior cover $f: Y \to X$ can be constructed as follows.  Let $N_{0} \subset N$ be the sublattice generated by $\{ v_{\rho} \}$.  We define $N'$ to be the prime-to-$p$ saturation, that is,
    \begin{equation*}
        N' = \{ v \in N \, | \, mv \in N_{0} \textrm{ for some }m \textrm{ with }\gcd(m,p)=1\}.
    \end{equation*}
    Let $\Sigma'$ be the fan given by intersecting $\Sigma$ with $N'$ and let $Y$ be the normal projective toric variety defined by $\Sigma' \subset N'$.  The map $\Sigma' \to \Sigma$ induced by $N' \subset N$ defines a purely inseparable finite morphism $f: Y \to X$.  Since $N'/N_{0}$ has order prime-to-$p$, the map $N_{0} \otimes \mathbb{F}_{p} \to N' \otimes \mathbb{F}_{p}$ is an isomorphism.  As the lattice $N_{0}$ is generated by the primitive generators of rays of $\Sigma'$,  we see that $N' \otimes \mathbb{F}_{p}$ satisfies Proposition \ref{prop:toricexam}.(3) so that $T_{Y}$ has the PQP.  Since $\Sigma' \to \Sigma$ preserves all the ray generators, it identifies $K_{Y} = f^{*}K_{X}$ so that $f$ is crepant.  In combination these properties show that $f$ is the superior cover by Theorem \ref{theo:pqpcanonicalcover}.  %\brian{AI gives a toric example where we cannot get the superior cover after a single height $1$ cover, but it requires multiple foliation steps.  The height $k$ cover is given by looking at $N_{0} + p^{k}N$.  For $k$ sufficiently large this will stabilize to $N'$, particularly, it will do this once $p^{k}$ kills the quotient.}
\end{construction}

\subsection{Koll\'ar's examples} 
\cite{Kollar95} gives some additional examples of projective globally $F$-regular varieties $X$ such that $T_{X}$ does not have the PQP.  (More precisely, such examples are the ``boundary'' case of Koll\'ar's construction.)  We verify that in these examples $Q_{max}$ is a torsion-free sheaf of rank $\dim X-1$ that is strongly stable (in the sense of \cite[Definition 5.3]{JLR26}) and is not numerically flat. 

We used Claude Opus 5.5 LLM for the computations in the following example.

\begin{example} \label{exam:kollar}
    We work over an algebraically closed field of odd characteristic $p$.  Choose a positive integer $n \geq 3$ so that $n+1=pd$ for some positive integer $d$.  Fix a general equation $f \in H^{0}(\mathbb{P}^{n},\mathcal{O}_{\mathbb{P}^{n}}(pd))$ and let $g: X \to \mathbb{P}^{n}$ be the purely inseparable cover of $\mathbb{P}^{n}$ defined by the equation $y^{p} = f(x_{0},\ldots,x_{n})$ in $\mathbb{P}(1^{n+1},d)$.  \cite{Kollar95} verifies that $X$ is a Fano variety with terminal Gorenstein singularities and that there is an exact sequence
    \begin{equation*}
        0 \to g^{*}\mathcal{O}_{\mathbb{P}^{n}}(-pd) \xrightarrow{\cdot df} g^{*}\Omega_{\mathbb{P}^{n}} \to \Omega_{X} \to g^{*}\mathcal{O}_{\mathbb{P}^{n}}(-d) \to 0.
    \end{equation*}
    Thus $-K_{X} = d g^{*}H$ where $H$ denotes the hyperplane class.  Dualizing, we get an exact sequence
    \begin{equation*}
        0 \to g^{*}\mathcal{O}_{\mathbb{P}^{n}}(d) \to T_{X} \to g^{*}T_{\mathbb{P}^{n}} \xrightarrow{\langle df, \cdot \rangle} g^{*}\mathcal{O}_{\mathbb{P}^{n}}(pd)
    \end{equation*}
    where the last map is surjective away from a finite set of points (due to the generality of $f$).

    We first verify that $X$ is globally $F$-regular when $f$ is general. Note that $X \subset \mathbb{P}(1^{n+1},d)$ avoids the singularities of the weighted projective space so that $\mathcal{O}(1)$ is an ample line bundle.  By \cite[Proposition 5.3]{SS10} it suffices to show that the homogeneous coordinate ring $R = \mathbb{K}[x_{0},\ldots,x_{n},y]/(y^{p}-f)$ is strongly $F$-regular.  Choose an index $j$ such that $c := \frac{\partial f}{\partial x_{j}}$ is non-zero.  Then $R_{c}$ is regular.  By \cite[Theorem 3.1]{Glassbrenner96}, it suffices to check that $c(y^{p}-f)^{q-1} \not \in \mathfrak{m}^{[q]}$ for some $q = p^{e}$ where $\mathfrak{m} \subset \mathbb{K}[x_{0},\ldots,x_{n},y]$ is the maximal ideal.  We will verify this when $p^{e-1} > \deg(c)$.  We can write
    \begin{align*}
        c(y^{p}-f)^{q-1} & =  c \sum_{i=0}^{q-1} (-1)^{q-1-i} {q-1 \choose i} y^{pi} f^{q-1-i}
    \end{align*}
    As $\mathfrak{m}^{[q]}$ is monomial, it suffices to find a single term that is not contained in $\mathfrak{m}^{[q]}$.
    Consider the term $\pm {q-1 \choose p^{e-1}-1} c y^{q-p} f^{q-p^{e-1}}$.  The coefficient is non-zero by Lucas's theorem.  By generality, we may assume that the coefficient of $(x_{0}\ldots x_{n})^{p-1}$ in $f^{p-1}$ is not zero.  Since $p^{e-1} > \deg(c)$, we have that $cy^{q-p}f^{q-p^{e-1}} = cy^{q-p}(f^{p-1})^{p^{e-1}}$ does not lie in $\mathfrak{m}^{[q]}$, finishing the verification.

    Let $Q$ be the cokernel of the map $g^{*}\mathcal{O}_{\mathbb{P}^{n}}(d) \to T_{X}$.  It is clear that $Q$ has numerically trivial first Chern class; since the rank of $Q$ is one less than the rank of $T_{X}$, we must have $Q = Q_{max}$.  We next verify that $Q$ is not numerically flat.   Indeed, let $F$ be the kernel of the map $T_{\mathbb{P}^{n}} \to \mathcal{O}_{\mathbb{P}^{n}}(pd)$.  Since this map is surjective away from a codimension $\geq 3$ closed subset, an easy computation shows that $c_{2}(F) = \frac{pd(pd-1)}{2}H^{2}$.  Since $g$ is faithfully flat, we have $Q \cong g^{*}F$ so that $Q$ is also not numerically flat.

    Finally, we verify that $Q$ is strongly stable.  Since $Q \cong g^{*}F$, \cite[Lemma 5.4]{JLR26} shows that it suffices to prove that $F$ is strongly stable. Note that a destabilizing subsheaf of $F_{abs}^{e\widehat{*}}F$ would have first Chern class of non-negative degree.  Let $U \subset \mathbb{P}^{n}$ be the big open subset where $F$ is locally free.  Then it suffices to prove that $H^{0}(U, \Lambda^{j} F_{abs}^{e*}F|_{U}) = 0$ for all $1 \leq j \leq n-2$ and all $e \geq 0$.  Since $F^{e}_{abs,*}\mathcal{O}_{\mathbb{P}^{n}}$ is a direct sum of line bundles with exactly one copy of $\mathcal{O}_{\mathbb{P}^{n}}$ and no other globally generated summands, by the projection formula this is equivalent to showing $H^{0}(U, \Lambda^{j} F|_{U}) = 0$.  Finally, as the first Chern class of $F$ vanishes, we have $\Lambda^{j} F|_{U} \cong \Lambda^{n-1-j} F^{\vee}|_{U}$ and it suffices to show the latter sheaves have no sections.

    Recall that $F^{\vee}$ is the reflexive hull of the cokernel of the map $\mathcal{O}_{\mathbb{P}^{n}}(-pd) \xrightarrow{df} \Omega_{\mathbb{P}^{n}}$.  Thus for $1 \leq j \leq n-2$ the Koszul complex for $df$ gives an exact sequence
    \begin{equation*}
        0 \to \mathcal{O}_{U}(-jpd) \to \Omega_{U}(-(j-1)pd) \to \ldots \to \Omega_{U}^{j-1}(-pd) \to \Omega_{U}^{j} \to \Lambda^{j} F^{\vee} |_{U} \to 0.
    \end{equation*}
    Since $\mathbb{P}^{n} \backslash U$ is finite, we have $H^{i}(U,\Omega_{U}^{j-i}(-ipd)) = H^{i}(\mathbb{P}^{n},\Omega_{\mathbb{P}^{n}}^{j-i}(-ipd)) = 0$ for $1 \leq j \leq n-2$ and $0 \leq i \leq j$. Thus, we see that $\Lambda^{j} F^{\vee} |_{U}$ has no global sections as desired.
\end{example}

\bibliographystyle{alpha}
\bibliography{freecurvescharp}
\end{document}